\documentclass[a4paper,12pt]{article}
\usepackage{amsmath}
\usepackage{amsfonts}
\usepackage{amssymb}
\usepackage{latexsym}
\usepackage{epsfig}
\usepackage{graphicx}
\usepackage{oldgerm}
\usepackage{theorem}

\def\C{\mathbb{C}}
\def\R{\mathbb{R}}
\def\N{\mathbb{N}}
\def\W{\mathbb{W}}
\def\Z{\mathbb{Z}}
\def\I{\mathbb{I}}
\def\J{\mathbb{J}}

\def\P{{\mathbb P}}
\def\E{{\mathbb E}}

\def\mbK{\mathbb{K}}

\def\bK{{\bf K}}

\def\bE{{\bf E}}

\def\rP{{\rm P}}
\def\rE{{\rm E}}
\def\rC{{\rm C}}
\def\rS{{\rm S}}

\def\1{{\bf 1}}
\def\0{{\bf 0}}
\def\x{\mib{x}}
\def\y{\mib{y}}
\def\u{\mib{u}}
\def\v{\mib{v}}
\def\w{\mib{w}}
\def\f{\mib{f}}
\def\t{\mib{t}}
\def\vchi{\mib{\chi}}

\def\X{\mib{X}}
\def\Y{\mib{Y}}
\def\bW{\mib{W}}
\def\bM{{\bf M}}
\def\bJ{{\bf J}}

\def\bA{{\bf A}}

\def\mM{\mathfrak{M}}

\def\supp{{\rm supp}\,}
\def\cG{{\cal G}}
\def\cM{{\cal M}}
\def\cF{{\cal F}}
\def\cI{{\cal I}}
\def\cS{{\cal S}}
\def\cD{{\cal D}}

\def\sc{{\sf c}}
\def\cC{{\cal C}}

\def\uj{\underline{j}}
\def\oj{\overline{j}}

\def\cA{{\cal A}}

\theorembodyfont{\itshape}

\newtheorem{thm}{Theorem}[section]
\newtheorem{lem}[thm]{Lemma}
\newtheorem{cor}[thm]{Corollary}
\newtheorem{prop}[thm]{Proposition}

\newtheorem{df}[thm]{Definition}

\newcommand{\mib}[1]{\mbox{\boldmath $#1$}}
\newcommand{\SSC}[1]{\section{#1}\setcounter{equation}{0}}
\newcommand{\qed}{\hbox{\rule[-2pt]{3pt}{6pt}}}

\begin{document}

\title{\bf 
Determinantal Martingales \\
and Noncolliding Diffusion Processes}
\author{
Makoto Katori
\footnote{
Department of Physics,
Faculty of Science and Engineering,
Chuo University, 
Kasuga, Bunkyo-ku, Tokyo 112-8551, Japan;
e-mail: katori@phys.chuo-u.ac.jp
}}
\date{9 July 2014}
\pagestyle{plain}
\maketitle
\begin{abstract}
Two aspects of noncolliding diffusion processes have been
extensively studied. One of them is the fact that
they are realized as harmonic Doob transforms of 
absorbing particle systems in the Weyl chambers.
Another aspect is integrability in the sense that
any spatio-temporal correlation function 
can be expressed by a determinant.
The purpose of the present paper is to clarify the connection
between these two aspects.
We introduce a notion of determinantal martingale
and prove that, if the system has
determinantal-martingale representation,
then it is determinantal.
In order to demonstrate the direct connection between the two aspects, 
we study three processes.
\end{abstract}



\normalsize

\SSC{Introduction \label{sec:Introduction}}

Two aspects of noncolliding diffusion processes have been
extensively studied in probability theory and random matrix theory.
One of them is the fact that,
although they are originally introduced as 
eigenvalue processes of matrix-valued diffusions \cite{Dys62,Bru91}, 
they are realized as harmonic Doob transforms of 
absorbing particle systems in the Weyl chambers \cite{Gra99,KO01,KT04}.
It implies that the noncolliding particle systems can be
regarded as multivariate extensions of the three-dimensional Bessel process,
which is realized as an $h$-transform of absorbing Brownian motion
in one-dimension with an absorbing wall at the origin \cite{KT11b}.
Another aspect is integrability in the sense that
any spatio-temporal correlation function can be explicitly expressed by a determinant
specified by a single continuous function called the correlation kernel.
Such determinantal processes are considered to be dynamical extensions
of random matrix models in which the eigenvalue ensembles provide
determinantal point processes \cite{Koe05,For10,Tao12}.
The purpose of the present paper is to clarify the connection
between these two aspects.
We introduce a notion of {\it determinantal martingale}
and prove that, if the interacting particle system has
an expression called {\it determinantal-martingale representation (DMR)},
then it is determinantal.
In order to demonstrate the direct connection between the two aspects
through DMR, we study three processes,
the noncolliding Brownian motion 
(Dyson's Brownian motion model with $\beta=2$),
the noncolliding squared Bessel process,
and the noncolliding Brownian motion on a circle.

We consider a continuous time Markov process $Y(t), t \in [0, \infty)$ 
on the state space $S$
which is a connected open set in $\R$.
It is a diffusion process in $S$ or a process showing
a position on the circumference 
$S=[0, 2 \pi r)$ of a diffusion process moving around
on the circle with a radius $r>0$;
$\rS^{1}(r) \equiv \{x \in [0, 2 \pi r) : x+2 \pi r = x\}$.
The probability space is denoted by $(\Omega, \cF, \rP_v)$
with expectation $\rE_v$, 
when the initial state is fixed to be $v \in S$.
When $v$ is the origin, the subscript is omitted.
We introduce a filtration $\{\cF(t) : t \in [0, \infty) \}$
generated by $Y$ so that it satisfies the usual conditions
(see, for instance, p.45 in  \cite{RY05}).
We assume that the process has
a transition density (TD),
$p(t, y|x)$, 
$t \in [0, \infty), x, y \in S$
such that for any measurable bounded 
function $f(t, x), t \in [0, \infty), x \in S$,
\begin{equation}
\rE[f(t, Y(t)) | \cF(s)]
= \int_{S} dy \, f(t, y) p(t-s, y|Y(s)) \quad
\mbox{a.s., \, $0 \leq s \leq t < \infty$}.
\label{eqn:def_p}
\end{equation}

For $N \in \N \equiv \{1,2, \dots\}$, we put
$$
\W_N(S)=\{\x=(x_1, \dots, x_N) \in S^{N} : x_1 < x_2 < \cdots < x_N \}.
$$
If $S=\R$, $\W_N(\R)$ is the {\it Weyl chamber} of type $A_{N-1}$ \cite{Gra99}.
For $\u=(u_1, \dots, u_N) \in \W_N(S)$,
we define a measure $\xi$ by a sum of
point masses concentrated on $u_j, 1 \leq j \leq N$; 
$
\xi(\cdot)=\sum_{j=1}^N \delta_{u_j}(\cdot).
$
Depending on $\xi$, we 
assume that there is a one-parameter
family of continuous functions
$$
\cM_{\xi}^v(\cdot, \cdot) : 
[0, \infty) \times S \mapsto \R
$$
with parameter $v \in S$, 
such that the following conditions are satisfied.
\begin{description}
\item{(M1)} \quad
$\cM_{\xi}^{u_k}(t, Y(t))$, $1 \leq k \leq N$, $t \in [0, \infty)$ are continuous-time martingales; 
$$
\rE[\cM_{\xi}^{u_k}(t, Y(t)) | \cF(s)]=\cM_{\xi}^{u_k}(s, Y(s)) \quad 
\mbox{a.s. for all $0 \leq s \leq t$.}
$$
\item{(M2)} \quad
For any time $t \geq 0$, $\cM_{\xi}^{u_k}(t, x), 1 \leq k \leq N$
are linearly independent functions of $x$.
\item{(M3)} \quad
For $1 \leq j, k \leq N$,
$$
\lim_{t \downarrow 0} \rE_{u_j}[\cM_{\xi}^{u_k}(t, Y(t))] 
= \delta_{jk}.
$$
\end{description}
We call $\cM_{\xi}(t, x)=\{\cM^v_{\xi}(t, x)\}$ martingale functions.

Let $\{Y_j(t), t \in [0, \infty) : 1 \leq j \leq N \}$
be a collection of $N$ independent copies 
of $Y(t), t \in [0, \infty)$.
We consider the $N$-component vector-valued Markov process
$\Y(t)=(Y_1(t), \dots, Y_N(t))$, $t \in [0, \infty)$,
for which the initial values are fixed to be
$Y_j(0)=u_j \in S, 1 \leq j \leq N$,
provided $\u=(u_1, \dots, u_N) \in \W_N(S)$, 
and
the probability space is denoted by
$(\Omega, \cF, \rP_{\u})$ with expectation
$\rE_{\u}$. 
We consider a determinant of the martingales
\begin{equation}
\cD_{\xi}(t, \Y(t)) 
=\det_{1 \leq j,k \leq N} \Big[ 
\cM_{\xi}^{u_k}(t, Y_j(t)) \Big],
\quad t \in [0, \infty). 
\label{eqn:D1}
\end{equation}
The condition (M2) is necessary so that 
it is not zero constantly.
This determinant is a continuous-time martingale 
and we call it a {\it determinantal martingale}.

Let $\mM(S)$ be the space of nonnegative integer-valued 
Radon measures on $S$,
which is a Polish space with the vague topology.
In this paper the cardinality of a finite set $A$ is denoted by $\sharp A$.
Any element $\xi$ of $\mM(S)$ can be represented as
$\xi(\cdot) = \sum_{j \in \I}\delta_{x_j}(\cdot)$
with a countable index set $\I$, in which
a sequence of points in $S$, $\x =(x_j)_{j \in \I}$, 
satisfies $\xi(K)=\sharp\{x_j, x_j \in K \} < \infty$ 
for any compact subset $K \subset S$.
We consider interacting particle systems 
as $\mM(S)$-valued processes and write them as
\begin{equation}
\Xi(t, \cdot)=\sum_{j=1}^N \delta_{X_j(t)}(\cdot),
\quad t \in [0, \infty),
\label{eqn:Xi1}
\end{equation}
where $\X(t)=(X_1(t), \dots, X_N(t)), t \in [0, \infty)$ 
are defined by a solution of
a given system of stochastic differential equations (SDEs).
We call $\x \in S^N$ a labeled configuration
and $\xi \in \mM(S)$ an unlabeled configuration.
The probability law of $\Xi(t, \cdot), t \in [0, \infty)$
starting from a fixed configuration $\xi \in \mM(S)$
is denoted by $\P_{\xi}$ 
and the process
specified by the initial configuration
is expressed by
$(\Xi, \P_{\xi})$.
The expectations with respect to $\P_{\xi}$
is denoted by $\E_{\xi}$.
We introduce a filtration $\{\cF_{\Xi}(t) : t \in [0, \infty)\}$
generated by $\Xi(t)$ satisfying the usual conditions.
Let $\rC_0(S)$ be the set of all continuous
real-valued functions with compact supports on $S$.
We set 
$$
\mM_{0}(S)= \{ \xi \in \mM(S) : 
\xi(\{x\}) \leq 1 \mbox { for any }  x \in S \},
$$
which gives a collection of configurations
without any multiple points.

Let $0 \leq t \leq T < \infty$.
We consider the expectation
of an $\cF_{\Xi}(t)$-measurable bounded function $F$,
$\E_{\xi} [F(\Xi(\cdot))]$.
It is sufficient 
to consider the case that $F$ is given as
\begin{equation}
F\left(\Xi(\cdot)\right)
= \prod_{m=1}^M g_m(\X(t_m))
\label{eqn:F_gm}
\end{equation}
for an arbitrary $M \in \N$, $0 \leq t_1< \cdots <t_M \leq T < \infty$
with bounded measurable functions $g_m$ 
on $S^N$, $1 \leq m \leq M$.
Since the particles are unlabeled 
in the process $(\Xi, \P_{\xi})$,
$g_m$'s are symmetric functions.

\begin{df}
\label{thm:DM_rep}
Given $\xi \in \mM_0(S)$, consider a process $(\Xi, \P_{\xi})$. 
If there exists a pair $(Y, \cM_{\xi})$
defining $\cD_{\xi}$ by (\ref{eqn:D1}) such that
for any $\cF_{\Xi}(t)$-measurable bounded function $F$,
$0 \leq t \leq T < \infty$, 
the equality 
$$
\E_{\xi}[F(\Xi(\cdot))]
=\rE_{\u} \left[
F \left( \sum_{j=1}^N \delta_{Y_j(\cdot)} \right) 
\cD_{\xi}(T, \Y(T)) \right]
$$
holds, 
then we say $(\Xi, \P_{\xi})$ has determinantal-martingale
representation (DMR)
associated with $(Y, \cM_{\xi})$.
\end{df}

For any integer $M \in \N$,
a sequence of times
$\t=(t_1,\dots,t_M) \in [0, \infty)^M$ with 
$0 \leq t_1 < \cdots < t_M < \infty$, 
and a sequence of functions
$\f=(f_{t_1},\dots,f_{t_M}) \in \rC_0(S)^M$,
the {\it moment generating function} of multi-time distribution
of $(\Xi, \P_{\xi})$ is defined by
\begin{equation}
\Psi_{\xi}^{\t}[\f]
\equiv \E_{\xi} \left[ \exp \left\{ \sum_{m=1}^{M} 
\int_{S} f_{t_m}(x) \Xi(t_m, dx) \right\} \right].
\label{eqn:GF1}
\end{equation}
It is expanded with respect to `test functions' 
$\chi_{t_m}(\cdot) \equiv e^{f_{t_m}(\cdot)}-1,
1 \leq m \leq M$ as
$$
\Psi_{\xi}^{\t}[\f]
=\sum_
{\substack
{0 \leq N_m \leq N, \\ 1 \leq m \leq M} }
\int_{\prod_{m=1}^{M} \W_{N_{m}}(S)}
\prod_{m=1}^{M} 
d \x_{N_m}^{(m)}
\prod_{j=1}^{N_{m}} 
\chi_{t_m} \Big(x_{j}^{(m)} \Big) 
\rho_{\xi} 
\Big( t_{1}, \x^{(1)}_{N_1}; \dots ; t_{M}, \x^{(M)}_{N_M} \Big),
$$
where $\x^{(m)}_{N_m}$ denotes
$(x^{(m)}_1, \dots, x^{(m)}_{N_m})$, 
and this expansion formula of $\Psi_{\xi}^{\t}[\f]$
defines the
{\it spatio-temporal correlation functions}
$\rho_{\xi}$ for the process $(\Xi, \P_{\xi})$.
Given an integral kernel, 
$\bK(s,x;t,y)$, 
$(s,x), (t,y) \in [0, \infty) \times S$,
a {\it Fredholm determinant} is defined as
\begin{eqnarray}
&& \mathop{{\rm Det}}_
{\substack{
(s,t)\in \{t_1, \dots, t_M\}^2, \\
(x,y)\in S^2}
}
 \Big[\delta_{st} \delta_x(\{y\})
+ \bK(s,x;t,y) \chi_{t}(y) \Big]
\nonumber\\
&& 
=\sum_
{\substack
{0 \leq N_m \leq N, \\ 1 \leq m \leq M} }
\int_{\prod_{m=1}^{M} \W_{N_{m}}(S)}
\prod_{m=1}^{M} 
d \x_{N_m}^{(m)}
\prod_{j=1}^{N_{m}} 
\chi_{t_m} \Big(x_{j}^{(m)} \Big) 
\det_{\substack
{1 \leq j \leq N_{m}, 1 \leq k \leq N_{n}, \\
1 \leq m, n \leq M}
}
\Bigg[
\bK(t_m, x_{j}^{(m)}; t_n, x_{k}^{(n)} )
\Bigg]
\nonumber\\
\label{eqn:Fredholm1}
\end{eqnarray}
for a sequence of functions
$\vchi=(\chi_{t_1},\dots,\chi_{t_M}) \in \rC_0(S)^M$, 
where 
$d \x^{(m)}_{N_m}= \prod_{j=1}^{N_m} dx^{(m)}_j$,
$1 \leq m \leq M$.
\begin{df}
\label{thm:determinantal}
If any moment generating function (\ref{eqn:GF1}) 
is given by a Fredholm determinant,
the process $(\Xi, \P_{\xi})$ is said to be determinantal. 
In this case, 
all spatio-temporal correlation functions
are given by determinants as 
\begin{equation}
\rho_{\xi} \Big(t_1,\x^{(1)}_{N_1}; \dots;t_M,\x^{(M)}_{N_M} \Big) 
=\det_{\substack
{1 \leq j \leq N_{m}, 1 \leq k \leq N_{n}, \\
1 \leq m, n \leq M}
}
\Bigg[
\mbK_{\xi}(t_m, x_{j}^{(m)}; t_n, x_{k}^{(n)} )
\Bigg],
\label{eqn:rho1}
\end{equation}
$0 \leq t_1 < \cdots < t_M < \infty$, 
$1 \leq N_m \leq N$,
$\x^{(m)}_{N_m} \in S^{N_m}, 1 \leq m \leq M \in \N$.
Here the integral kernel, $\mbK_{\xi}:([0, \infty) \times S)^2 \mapsto \R$,
is a function of initial configuration
$\xi$ and is called
the correlation kernel.
\end{df}

Let $\1(\omega)$ be the indicator of $\omega$;
$\1(\omega)=1$ if $\omega$ is satisfied,
and $\1(\omega)=0$ otherwise.
The main theorem of the present paper is the following.
\begin{thm}
\label{thm:DM_det}
If $(\Xi, \P_{\xi})$ has DMR associated with $(Y, \cM_{\xi})$, 
then it is determinantal
with the correlation kernel
\begin{equation}
\mbK_{\xi}(s,x;t,y)
=\int_{S} \xi(dv) p(s, x|v) \cM_{\xi}^{v}(t,y)
- \1(s>t) p(s-t,x|y),
\label{eqn:K1}
\end{equation}
$(s,x), (t,y) \in [0, \infty) \times S$,
where $p$ is the TD of the process $Y$.
\end{thm}

This type of asymmetric correlation kernel (\ref{eqn:K1}) 
was first obtained
by Eynard and Mehta 
for a multi-matrix model \cite{EM98}
and by Nagao and Forrester 
for the noncolliding Brownian motion started at a special
initial distribution (the GUE eigenvalue distribution) \cite{NF98}, 
and has been extensively studied
\cite{FNH99,TW04,Meh04,BR05,For10,KT10,KT13}.
Note that, while in determinantal point processes
correlation kernels are usually symmetric
\cite{Sos00,ST03,BKPV09}, 
the present `dynamical correlation kernel'  
is asymmetric with respect to
the exchange of two points $(s, x)$ and $(t, y)$
on the spatio-temporal plane $[0, \infty) \times S$ 
and shows causality in the system.

In order to show applications of Theorem \ref{thm:DM_det},
we consider the following three kinds of interacting particle systems,
$(\Xi, \P_{\xi}), \xi=\sum_{j=1}^N \delta_{u_j} \in \mM_0(S)$ with
$\Xi(t, \cdot)=\sum_{j=1}^N \delta_{X_j(t)}(\cdot), t \in [0, \infty)$.
For each system of SDEs, $B_j(t), 1 \leq j \leq N, t \geq 0$
denote a set of independent one-dimensional standard Brownian motions 
(BMs) started at 0.
(From now on, BM means a one-dimensional standard
Brownian motion unless specially mentioned.)  
\begin{description}
\item[Process 1 :] \,
Noncolliding Brownian motion (Dyson's Brownian motion model with $\beta=2$);
\begin{eqnarray}
S &=& \R,
\nonumber\\
X_j(t) &=& u_j+ B_j(t)+ 
\sum_{\substack{ 1 \leq k \leq N,\\ k \not= j}} \int_0^t
\frac{ds}{X_j(s)-X_k(s)},\quad 1 \leq j \leq N, \quad t \geq 0.
\label{eqn:noncollBM}
\end{eqnarray}

\item[Process 2 :] \,
Noncolliding squared Bessel process (BESQ$^{(\nu)}$) with $\nu >-1$; 
\begin{eqnarray}
S &=& \R_+ \equiv \{x \in \R: x \geq 0\},
\nonumber\\
X_j(t) &=& u_j+ \int_0^t 2 \sqrt{X_j(s)} 
d B_j(s) + 2 (\nu+1) t
\nonumber\\
&& 
+ \sum_{\substack{1 \leq k \leq N, \\ k \not= j}}
\int_0^t
\frac{4 X_j(s)ds}{X_j(s)-X_k(s)},
\quad 1 \leq j \leq N, \quad t \geq 0,
\label{eqn:noncollBESQ}
\end{eqnarray}
where if $-1 < \nu < 0$
the reflection boundary condition is assumed at
the origin.

\item[Process 3 :] \,
Noncolliding BM on a circle with a radius $r>0$;
(trigonometric extension of Dyson's Brownian motion model
with $\beta=2$); \\
We solve the SDEs
\begin{equation}
\check{X}_j(t) = u_j+ B_j(t)
+\frac{1}{2 r} \sum_{\substack{1 \leq k \leq N, \cr k \not= j}}
\int_0^t
\cot \left( \frac{\check{X}_j(s)-\check{X}_k(s)}{2 r} \right) ds,
\quad 1 \leq j \leq N, \quad t \geq 0,
\label{eqn:noncollBM_S1}
\end{equation}
on $\R$ and then define the process on
$S=[0, 2 \pi r)$
by 
\begin{equation}
X_j(t)=\check{X}_j(t) \quad
\mbox{mod $2 \pi r$},
\quad 1 \leq j \leq N, \quad t \geq 0.
\label{eqn:noncollBM_S1b}
\end{equation}

Note that $\cot (x/2r)$ is a periodic function of $x$ with
period $2 \pi r$.
By definition (\ref{eqn:noncollBM_S1b}),
measurable functions for Process 3 should be
periodic with period $2 \pi r$ in the following sense.
For $0 \leq t < \infty$, if an $\cF_{\Xi}(t)$-measurable
function $F$ is given in the form (\ref{eqn:F_gm}),
then for any $n \in \Z$,
$$
g_m((x_j+2 \pi r n)_{j=1}^N)
=g_{m}(\x), 
\quad 1 \leq m \leq M.
$$
\end{description}
Provided $\u=\X(0) \in \W_N(S)$, let
$\tau_{\rm collision}= \inf \{ t > 0: \X(t) \notin \W_N(S)\}$.
It has been proved that
$\P_{\xi}(\tau_{\rm collision}=\infty)=1$
for these three processes
\cite{Cha92,RS93,GM13,CL01}.
The noncolliding Brownian motion was introduced by
Dyson as the eigenvalue process of Hermitian
matrix-valued process \cite{Dys62}
and has been also called Dyson's Brownian motion model
(with parameter $\beta=2$) 
\cite{Spo87,NF98,Gra99,Joh01,Meh04,TW04,KT04,KT09,KT10,Osa12,KT13,Osa13,Osa13b}.
The noncolliding BESQ$^{(\nu)}$ was studied by
K\"onig and O'Connell \cite{KO01},
where it is considered as the generalization of
the eigenvalue process of a matrix-valued process
called the Laguerre process (or the complex Wishart process).
See also \cite{Bru91,Gra99,FNH99,KT04,TW07,KT07a,KT11}.
The system (\ref{eqn:noncollBM_S1}) 
with the identification (\ref{eqn:noncollBM_S1b}) can be regarded as
a dynamical extension of the circular unitary ensemble (CUE)
of random matrix theory (see Section 11.8 in \cite{Meh04}
and Chapter 11 in \cite{For10}) \cite{HW96,CL01}.
The dynamics was studied in \cite{NF02}
and papers cited therein. 
It is interesting to see that
(a one-parameter extension of) (\ref{eqn:noncollBM_S1})
is discussed as a driving system
for a multiple radial Schramm-Loewner evolution
by Cardy \cite{Car03}.
Quite recently many interesting properties of 
related processes on a circle are reported in \cite{PSS10,PS11,SMCF13,LW13}.
Process 3 is a trigonometric extension of Process 1
and in the limit $r \to \infty$
Process 3 should be reduced to Process 1.
Since functions used to represent Process 3 in the present paper
are all analytic with respect to $r$,
the hyperbolic extension will be similarly discussed \cite{CL01}.

We introduce sets of entire functions of $z \in \C$ \cite{Lev96}, 
$1 \leq k \leq N$, for
$\u \in \W_N(S)$, $\xi=\sum_{j=1}^N \delta_{u_j}$, $v \in S$,
$r>0$, 
\begin{equation}
\Phi_{\xi}^{v}(z)
=\left\{ \begin{array}{ll}
\displaystyle{
\prod_{\substack{1 \leq \ell \leq N, \cr u_{\ell} \not=v}}
\frac{z-u_{\ell}}{v-u_{\ell}},
}
& \quad \mbox{for Processes 1 and 2}, \cr
\displaystyle{
\prod_{\substack{1 \leq \ell \leq N, \cr u_{\ell} \not=v}}
\frac{\sin((z-u_{\ell})/2r)}{\sin((v-u_{\ell})/2r)},
}
& \quad \mbox{\rm for Process 3}, \cr
\end{array} \right.
\label{eqn:Phi}
\end{equation}
and integral transformations of function $f=f(W)$,  
\begin{equation}
\cI[f(W)|(t,x)]
= \left\{ \begin{array}{ll}
\displaystyle{
\int_{\R} dw \, f(i w) q(t,w|x),
}
& \quad \mbox{for Processes 1 and 3}, \cr
\displaystyle{
\int_{\R_+} dw \, f(-w) q^{(\nu)}(t,w|x),
}
& \quad \mbox{for Process 2},
\end{array} \right.
\label{eqn:I}
\end{equation}
with 
\begin{eqnarray}
\label{eqn:q}
q(t, w|x) &=&  
\frac{e^{-(ix+w)^2/2t}}{\sqrt{2 \pi t}}, 
\\
\label{eqn:q_nu}
q^{(\nu)}(t,w|x) &=& 
\left( \frac{w}{x} \right)^{\nu/2}
\frac{e^{(x-w)/2t}}{2t}
J_{\nu} \left( \frac{\sqrt{xw}}{t} \right),
\quad \nu > -1, 
\end{eqnarray}
where $i=\sqrt{-1}$ and
$J_{\nu}$ is the Bessel function
$$
J_{\nu}(z) = \sum_{n=0}^{\infty} 
\frac{(-1)^n}{\Gamma(n+1) \Gamma(n+1+\nu)}
\left( \frac{z}{2} \right)^{2n+\nu}
$$
with the Gamma function $\Gamma(z)$,
which is an analytic function of $z$ in a cut plane 
$\C \setminus (-\infty, 0)$ \cite{Wat44}.
Appropriate limits should be taken for $t=0$ in (\ref{eqn:q})
and (\ref{eqn:q_nu}) and for $x=0$ in (\ref{eqn:q_nu})
as explained in Section \ref{sec:martingales}.
(Note that $W$ in LHS of (\ref{eqn:I}) is a dummy variable,
but it will be useful in order to specify a function $f$,
especially when we extend the integral transformation for
multivariate functions in Section \ref{sec:martingales}.)

Corresponding to the three interacting
$N$-particle systems, 
we consider the following three kinds of
one-dimensional processes.
The first one is BM on $S=\R$, whose TD 
started at $x \in \R$ is given by
\begin{equation}
p(t, y|x)=
p_{\rm BM}(t, y|x) =  \left\{ \begin{array}{ll}
\displaystyle{
\frac{e^{-(y-x)^2/2t}}{\sqrt{2 \pi t}}},
& t>0, y \in \R \cr
& \cr
\delta_x(\{y\}),
& t=0, y \in \R.
\end{array} \right.
\label{eqn:p}
\end{equation}
The second one is BESQ$^{(\nu)}$ with $\nu > -1$
on $S=\R_+$, which is 
given by the solution of the SDE, 
$$
Y(t)=u+\int_0^t 2 \sqrt{Y(s)} d B(s)
+2(\nu+1) t, \quad t \geq 0, \quad u > 0, 
$$
where $B$ is BM, and 
if $-1 < \nu < 0$ a reflecting wall
is put at the origin.
The TD is given by
\begin{eqnarray}
p(t,y|x) &=& p^{(\nu)}(t, y|x)
\nonumber\\
&=& \left\{ \begin{array}{ll}
\displaystyle{
\left( \frac{y}{x} \right)^{\nu/2}
\frac{e^{-(x+y)/2t}}{2t}
I_{\nu} \left( \frac{\sqrt{xy}}{t} \right)},
& t>0, x>0, y \in \R_+, \cr
\displaystyle{
\frac{y^{\nu}e^{-y/2t}}{(2t)^{\nu+1} \Gamma(\nu+1)}},
& t >0, x=0, y \in \R_+, \cr
& \cr
\delta_x(\{y\}),
& t=0, x, y \in \R_+,
\end{array} \right.
\label{eqn:p_nu}
\end{eqnarray}
where $I_{\nu}(x)$ is the modified Bessel function
of the first kind \cite{Wat44}
$$
I_{\nu}(x) = \sum_{n=0}^{\infty} 
\frac{(x/2)^{2n+\nu}}{\Gamma(n+1) \Gamma(n+1+\nu)}, \quad x \in \R_+.
$$
The third one is a Markov process on $S=[0, 2\pi r)$,
$r >0$, whose TD is given by
\begin{equation}
p(t, y|x)= p^r(t, y|x; N)= \left\{
\begin{array}{ll}
\displaystyle{\sum_{\ell \in \Z} p_{\rm BM}(t, y+2 \pi r \ell|x)}, \quad & \mbox{if $N$ is odd}, \cr
\displaystyle{\sum_{\ell \in \Z} (-1)^{\ell} 
p_{\rm BM}(t, y+2 \pi r \ell|x)}, \quad & \mbox{if $N$ is even},
\end{array} \right.
\label{eqn:prB1}
\end{equation}
$x, y \in [0, 2 \pi r)$, $t \geq 0$,
where $N$ is the number of particles of Process 3.
This process shows a position on the
circumference $S=[0, 2 \pi r)$ of a Brownian motion
moving around on the circle $\rS^1(r)$
(with alternating signed densities if $N$ is even).
We will explain in Section \ref{sec:alcove} that
it represents an elementary process
of the system of $N$ Brownian motions
on $\rS^1(r)$, which is killed when any collision
of particles occurs.
In the following, we call this one-dimensional Markov process
$Y(t) \in [0, 2 \pi r)$, $t \in [0, \infty)$
simply `BM on $[0, 2 \pi r)$'.

We will prove the following.
\begin{thm}
\label{thm:3_process}
For any $\xi \in \mM_0(S)$, 
the three processes have DMRs 
associated with $(Y, \cM_{\xi})$
such that
$$
\mbox{$Y$ is given by} \, \left\{ \begin{array}{ll}
\mbox{BM on $\R$}, & \mbox{for Process 1}, \cr
\mbox{BESQ$^{(\nu)}$ on $\R_+$}, & \mbox{for Process 2}, \cr
\mbox{BM on $[0, 2 \pi r)$}, & \mbox{for Process 3}, \cr
\end{array} \right. 
$$
and $\cM_{\xi}=\{\cM_{\xi}^v(t, x)\}, v \in S$ is given by
\begin{equation}
\cM^{v}_{\xi}(t, x)
=\cI[\Phi_{\xi}^{v}(W)|(t,x)].
\label{eqn:cMs}
\end{equation}
Then they are all determinantal 
with correlation kernels (\ref{eqn:K1}) with (\ref{eqn:cMs}).
\end{thm}

In the present paper, 
Section \ref{sec:proof1} is devoted to proof of Theorem \ref{thm:DM_det}.
In Section \ref{sec:martingales} properties of 
$(Y, \cM_{\xi})$ are studied for the three processes. 
Theorem \ref{thm:3_process} is proved in Section \ref{sec:proof2},
where direct relations between DMRs and $h$-transforms are clarified.
It is shown in Section \ref{sec:multi} that,
even though the systems of SDEs (\ref{eqn:noncollBM})-(\ref{eqn:noncollBM_S1})
cannot be solved if initial configurations have
multiple points, $\xi \in \mM(S) \setminus \mM_0(S)$,
the martingale functions $\cM_{\xi}$ can be 
rewritten in the form so that 
they are valid also for $\xi \in \mM(S) \setminus \mM_0(S)$.
The spatio-temporal correlations of determinantal processes
will determine the entrance laws for the processes
from $\xi \in \mM(S) \setminus \mM_0(S)$.
In Section \ref{sec:relax}, 
as an application of Theorem \ref{thm:3_process},
we show for Process 3 that a non-equilibrium dynamics
exhibiting a {\it relaxation phenomenon}
is well-described by the obtained correlation kernel.

Consider Processes 1 and 3. 
Since $\Phi_{\xi}^v(z)$ given by (\ref{eqn:Phi}) are
entire functions of $z$, 
their integral transformations (\ref{eqn:I}) 
with the integral kernel (\ref{eqn:q})
are rewritten as
\begin{equation}
\cI[\Phi_{\xi}^v(W)|(t, x)]
= \int_{\R} d \widetilde{y} \, 
\Phi_{\xi}^v(x+i \widetilde{y}) p_{\rm BM}(t, \widetilde{y}|0)
\label{eqn:cal1}
\end{equation}
by changing the integral variable as 
$w \to \widetilde{y}$ with
$\widetilde{y}=ix+w$
and by Cauchy's integral theorem.
Let $\widetilde{Y}(t), t \in [0, \infty)$ be
BM on $\R$ started at 0, which is independent from
$Y(t), t \in [0, \infty)$.
We write the expectation with respect to $\widetilde{Y}(\cdot)$
as $\widetilde{\rm E}$.
Define a complex process
$$
Z(t)=Y(t)+i \widetilde{Y}(t),
\quad t \in [0, \infty).
$$
For Process 1 it is a complex Brownian motion
on $\C$, and for Process 3 it is 
that defined on 
$\{z=x+i y : x \in [0, 2 \pi r), y \in \R \}, r >0$.
For Processes 1 and 3, (\ref{eqn:cal1}) implies the equality
$\cI[\Phi_{\xi}^v(W)|(t, Y(t))]
= \widetilde{\rE}[\Phi_{\xi}^v(Z(t))], 
t \in [0, \infty)$. 
Then (\ref{eqn:cMs}) gives
\begin{equation}
\cM_{\xi}^{v}(t, Y(t))
=\widetilde{\rm E}[\Phi_{\xi}^v(Z(t))], \quad 
t \in [0, \infty).
\label{eqn:CBM}
\end{equation}
In these cases, the present DMR is reduced to be
the {\it complex Brownian motion representation (CBMR)},
which was introduced in \cite{KT13} for Process 1.
On the other hand, Process 2
does not have CBMR.
We find, however, that if we transform
BESQ$^{(\nu)}$, $Y^{(\nu)}(\cdot)$, to
the Bessel process
(BES$^{(\nu)}$), $\widehat{Y}^{(\nu)}(\cdot)$,  by 
$\widehat{Y}^{(\nu)}=\sqrt{Y^{(\nu)}}$,
and if we set
$\nu=n+1/2, n \in \N_0 \cup \{-1 \}$,
where $\N_0 \equiv \N \cup \{0\}$,
then the corresponding martingales
$\widehat{\cM}_{\xi}$
can be represented by using a complex process
\begin{equation}
\widehat{Z}^{(n+1/2)}(t)=\widehat{Y}^{(n+1/2)}(t)+i \widetilde{Y}(t),
\quad t \in [0, \infty), 
\label{eqn:hatZ}
\end{equation}
where $\widetilde{Y}(t)$ is BM
defined independently from $\widehat{Y}^{(n+1/2)}(t)$, $t \in [0, \infty)$.
The {\it complex-process representation (CPR)} 
for the noncolliding BES$^{(n+1/2)}$, $n \in \N_0 \cup \{-1\}$
is reported in Section \ref{sec:CPR}.

\SSC{Proof of Theorem \ref{thm:DM_det} \label{sec:proof1}}
\subsection{Preliminaries \label{sec:proof1_pre}}
For $n \in \N$, an index set $\{1,2, \dots, n\}$
is denoted by $\I_{n}$.
Fixed $N \in \N$ with $N' \in \I_N$, 
we write $\J \subset \I_N, \sharp \J=N'$,
if $\J=\{j_1, \dots, j_{N'}\},
1 \leq j_1 < \dots < j_{N'} \leq N$.
For  $\x=(x_1, \dots, x_N) \in S^N$, 
put $\x_{\J}=(x_{j_1}, \dots, x_{j_{N'}})$.
In particular, we write
$\x_{N'}=\x_{\I_{N'}}, 1 \leq N' \leq N$.
(By definition $\x_N=\x$.)

For a finite index set $\J$, we write the collection of all permutations
of elements in $\J$ as $\cS(\J)$.
In particular, 
for $\I_n=\{1,2,\dots, n \}, n \in N$, 
we express $\cS(\I_{n})$ simply by $\cS_{n}$.
For $\x=(x_1, \dots, x_N) \in S^n$, $\sigma \in \cS_n$, 
we put $\sigma(\x)=(x_{\sigma(1)}, \dots, x_{\sigma(n)})$.
For an $n \times n$ matrix $B=(B_{jk})_{1 \leq j,k \leq n}$, 
the determinant is defined by
\begin{eqnarray}
\det B &=& \det_{1 \leq j, k \leq n} [B_{jk}]
\nonumber\\
&=& \sum_{\sigma \in \cS_{n}} {\rm sgn}(\sigma)
\prod_{j=1}^n B_{j \sigma(j)}.
\label{eqn:det1}
\end{eqnarray}
Any permutation
$\sigma \in \cS_n$ can be decomposed
into a product of cyclic permutations.
Let the number of cyclic permutations in the decomposition
be $\ell(\sigma)$ and express $\sigma$ by
$$
\sigma = \sc_1 \sc_2 \cdots \sc_{\ell(\sigma)},
$$
where $\sc_{\lambda}$
denotes a cyclic permutation
$$
\sc_{\lambda}= 
(c_\lambda(1) c_\lambda(2) \cdots c_\lambda(q_{\lambda}) ), \quad
1 \leq q_{\lambda} \leq n,  
\quad1 \leq \lambda \leq \ell(\sigma).
$$
For each $1 \leq \lambda \leq \ell(\sigma)$,
we write the set of entries 
$\{c_{\lambda}(j)\}_{j=1}^{q_{\lambda}}$ of $\sc_{\lambda}$
simply as $\{\sc_{\lambda}\}$,
in which the periodicity
$c_{\lambda}(j+q_{\lambda})=c_{\lambda}(j), 1 \leq j \leq q_{\lambda}$ 
is assumed.
By definition, for each $1 \leq \lambda \leq \ell(\sigma)$,
$c_{\lambda}(j), 1 \leq j \leq q_{\lambda}$
are distinct indices chosen from $\I_n$, and
$\{\sc_{\lambda}\} \cap \{\sc_{\lambda'}\} = \emptyset$
for $1 \leq \lambda \not= \lambda' \leq \ell(\sigma)$.
By definition 
$\sum_{\lambda=1}^{\ell(\sigma)} q_{\lambda}=n$.
The determinant (\ref{eqn:det1}) is also expressed as
\begin{equation}
\det B = \sum_{\sigma \in \cS_n}
(-1)^{n-\ell(\sigma)}
\prod_{\lambda=1}^{\ell(\sigma)}
\prod_{j=1}^{q_{\lambda}} 
B_{c_{\lambda}(j) c_{\lambda}(j+1)}.
\label{eqn:det2}
\end{equation}

\subsection{Lemmas \label{sec:proof1_lemma}}

\begin{lem}
\label{thm:reducibility}
Assume that
$\xi(\cdot)=\sum_{j=1}^N \delta_{u_j}(\cdot)$
with $\u \in \W_N(S)$.
Let $1 \leq N' \leq N$.
For $0 < t \leq T < \infty$ and an $\cF(t)$-measurable symmetric function
$F_{N'}$ on $S^{N'}$,
\begin{eqnarray}
&& \sum_{\J \subset \I_N, 
\sharp \J=N'}
\rE_{\u} \left[
F_{N'}(\Y_{\J}(t))
\cD_{\xi}(T, \Y(T)) \right]
\nonumber\\
&& \quad
= \int_{\W_{N'}(S)} \xi^{\otimes N'} (d\v)
\rE_{\v} \left[
F_{N'}(\Y_{N'}(t))
\cD_{\xi}(T, \Y_{N'}(T)) \right].
\label{eqn:reducibility}
\end{eqnarray}
\end{lem}

\begin{lem}
\label{thm:Fredholm}
Let $\u \in \W_N(S)$ and
$\xi=\sum_{j=1}^N \delta_{u_j}$.
Assume that there exists a pair $(Y, \cM_{\xi})$
satisfying conditions (M1)-(M3) and $\cD_{\xi}$
is defined by (\ref{eqn:D1}).
Then for any $M \in \N$, $0 \leq t_1 < \cdots < t_M \leq T < \infty$,
$\chi_{t_m} \in \rC_0(S), 1 \leq m \leq M$,
the equality
\begin{eqnarray}
&& \rE_{\u} \left[
\prod_{m=1}^M \prod_{j=1}^N
\{1+\chi_{t_m}(Y_j(t_m)) \}
\cD_{\xi}(T, \Y(T)) \right]
\nonumber\\
&& \qquad =
\mathop{{\rm Det}}_
{\substack{
(s,t)\in \{t_1, \dots, t_M\}^2, \\
(x,y)\in S^2}
}
 \Big[\delta_{st} \delta_x(\{y\})
+ \mbK_{\xi}(s,x;t,y) \chi_{t}(y) \Big]
\label{eqn:Fredholm2}
\end{eqnarray}
holds, 
where $\mbK_{\xi}$ is given by (\ref{eqn:K1}).
\end{lem}

Lemma \ref{thm:reducibility} shows the
{\it reducibility} of the determinantal martingale
in the sense that,
if we observe a symmetric function depending
on $N'$ variables, $N' \leq N$,
then the size of determinantal
martingale can be reduced from $N$ to $N'$.
The proof is given in Section \ref{sec:proof1_lemmaA}.

In Lemma \ref{thm:Fredholm}, the LHS of (\ref{eqn:Fredholm2})
is an expectation of a usual determinant multiplied by
test functions, while the RHS is a Fredholm determinant.
Note that the expectation in the LHS will be calculated
by performing integrals using TD, $p$ of the process $Y$
as an integral kernel,
while $p$ is involved in the integral representation
(\ref{eqn:K1}) of the correlation kernel $\mbK_{\xi}$
for the Fredholm determinant in the RHS.
Therefore, simply to say, this equality is just obtained
by applying Fubini's theorem.
Since the quantities in (\ref{eqn:Fredholm2}) are multivariate
and multi-time joint distribution is considered,
however, we also need combinatorics arguments
to prove Lemma \ref{thm:Fredholm}.
The proof is given in Section \ref{sec:proof1_lemmaB}.
It is quite similar to the proof for the CBMR
for Process 1 given in \cite{KT13}, but
some parts of the present proof are simpler,
since here we treat only real processes.

\subsection{Proof of Lemma \ref{thm:reducibility} \label{sec:proof1_lemmaA}}

By definition (\ref{eqn:D1}) with (\ref{eqn:det1}),
and independence of $Y_j(\cdot), 1 \leq j \leq N$,
the LHS of (\ref{eqn:reducibility}) is equal to
\begin{eqnarray}
&& \sum_{\J \subset \I_N, \sharp \J=N'}
\rE_{\u} \left[ F_{N'}(\Y_{\J}(t))
\det_{j, k \in \I_N}
[\cM_{\xi}^{u_k}(T, Y_j(T))] \right]
\nonumber\\
&& \, = \sum_{\J \subset \I_N, \sharp \J=N'}
\rE_{\u} \left[ F_{N'}(\Y_{\J}(t))
\sum_{\sigma \in \cS_N} {\rm sgn}(\sigma)
\prod_{j=1}^N \cM_{\xi}^{u_{\sigma(j)}}(T, Y_{j}(T)) \right]
\nonumber\\
&& \, = \sum_{\J \subset \I_N, \sharp \J=N'}
\sum_{\sigma \in \cS_N} {\rm sgn}(\sigma)
\rE_{\u} \left[ F_{N'}(\Y_{\J}(t))
\prod_{j \in \J} \cM_{\xi}^{u_{\sigma(j)}}(T, Y_j(T))
\prod_{k \in \I_N \setminus \J}
\cM_{\xi}^{u_{\sigma(k)}}(T, Y_k(T)) \right]
\nonumber\\
&& \, = \sum_{\J \subset \I_N, \sharp \J=N'}
\sum_{\sigma \in \cS_N} {\rm sgn}(\sigma)
\rE_{\u} \left[ F_{N'}(\Y_{\J}(t))
\prod_{j \in \J} \cM_{\xi}^{u_{\sigma(j)}}(T, Y_j(T)) \right]
\nonumber\\
&& \qquad \qquad \qquad \qquad \qquad \times
\prod_{k \in \I_N \setminus \J}
\rE_{\u} \left[ \cM_{\xi}^{u_{\sigma(k)}}(T, Y_k(T)) \right].
\label{eqn:A1}
\end{eqnarray}
By the condition (M1) of $\cM_{\xi}$, 
$$
\prod_{k \in \I_N \setminus \J}
\rE_{\u} \left[ \cM_{\xi}^{u_{\sigma(k)}}(T, Y_k(T)) \right]
=
\prod_{k \in \I_N \setminus \J}
\rE_{\u} \left[ \cM_{\xi}^{u_{\sigma(k)}}(t, Y_k(t)) \right],
\quad \forall t \in [0, T], 
$$
and the condition (M3) implies
it is equal to 
$\prod_{k \in \I_N \setminus \J}
\delta_{k \sigma(k)}$.
Then (\ref{eqn:A1}) becomes
\begin{eqnarray}
&&  \sum_{\J \subset \I_N, \sharp \J=N'}
\sum_{\sigma \in \cS(\J)} {\rm sgn}(\sigma)
\rE_{\u} \left[ F_{N'}(\Y_{\J}(t))
\prod_{j \in \J} \cM_{\xi}^{u_{\sigma(j)}}(T, Y_j(T)) \right]
\nonumber\\
&& \quad =
\sum_{\J \subset \I_N, \sharp \J=N'}
\rE_{\u} \left[ F_{N'}(\Y_{\J}(t))
\det_{j, k \in \J} [ \cM_{\xi}^{u_k}(T, Y_j(T))] \right]
\nonumber\\
&& \quad =
\int_{\W_{N'}(S)} \xi^{\otimes N'}(d \v)
\rE_{\v} \left[ F_{N'}(\Y_{N'}(t))
\det_{j, k \in \I_{N'}} [ \cM_{\xi}^{u_k}(T, Y_j(T))] \right],
\nonumber
\end{eqnarray}
where equivalence of $Y_j(t), t \in [0, \infty), 1 \leq j \leq N$ 
in probability law is used.
This is the RHS of (\ref{eqn:reducibility}) and 
the proof is completed. \qed

\subsection{Proof of Lemma \ref{thm:Fredholm} \label{sec:proof1_lemmaB}}

By performing binomial expansion of
$\prod_{m=1}^M \prod_{j=1}^N
\{1+\chi_{t_m}(Y_j(t_m)) \}$
and by Lemma \ref{thm:reducibility},
the LHS of (\ref{eqn:Fredholm2}) becomes
$$
\sum_
{\substack
{N_m \geq 0, \\ 1 \leq m \leq M} }
\sum_{1 \leq p \leq N}
\sum_{\substack{
\sharp \J_m=N_m, \\ 1 \leq m \leq M, \\
\bigcup_{m=1}^{M} \J_m= \I_{p}}}
\int_{\W_{p}(S)} \ 
\xi^{\otimes p}(d\v)
\rE_{\v} \left[
\prod_{m=1}^{M} 
\prod_{j_m\in \J_m} \chi_{t_m}(Y_{j_m}(t_m))
\cD_{\xi}(T, \Y_p(T))\Big] \right].
$$
On the other hand, the RHS of (\ref{eqn:Fredholm2})
has an expansion according to the definition
of Fredholm determinant (\ref{eqn:Fredholm1}).
Then, for proof of Lemma \ref{thm:Fredholm},
it is enough to show that, 
for any $M \in \N, (N_1, \dots, N_M) \in \N^M$, the equality
\begin{eqnarray}
&& \int_{\prod_{m=1}^{M} \W_{N_m}(S)} 
\prod_{m=1}^{M} \left\{ d \x_{N_m}^{(m)}
\prod_{j=1}^{N_m} \chi_{t_m} 
\Big(x_{j}^{(m)} \Big) \right\}
\det_{\substack{1 \leq j \leq N_{m}, 1 \leq k \leq N_{n},
\\ 1\le m, n \le M}}
\Bigg[
\mbK_{\xi}(t_m, x_{j}^{(m)}; t_{n}, x_{k}^{(n)} )
\Bigg]
\nonumber\\
&=&
\sum_{1 \leq p \leq N}
\sum_{\substack{
\sharp \J_m=N_m, \\ 1 \leq m \leq M, \\
\bigcup_{m=1}^{M} \J_m= \I_{p}}}
\int_{\W_{p}(S)} \ 
\xi^{\otimes p}(d\v)
\rE_{\v} \left[
\prod_{m=1}^{M} 
\prod_{j_m\in \J_m} \chi_{t_m}(Y_{j_m}(t_m))
\det_{j,k \in \I_p}
\Big[\cM_{\xi}^{v_k}(T, Y_j(T))\Big] \right]
\nonumber\\
\label{eqn:E=det}
\end{eqnarray}
is established. 
Here we will prove (\ref{eqn:E=det}) for
arbitrary chosen
$M \in \N$ and $(N_1, \dots, N_M) \in \N^M$.
The proof consists of four steps.

\vskip 0.3cm
\noindent
{\it Step 1: Restricted binomial expansion}
\vskip 0.3cm
Let $\I^{(1)}=\I_{N_1}$ and 
$\I^{(m)}=\I_{\sum_{j=1}^{m} N_{j}} 
\setminus \I_{\sum_{j=1}^{m-1} N_{j}}
=\{\sum_{j=1}^{m-1} N_j+1, \dots, \sum_{j=1}^m N_j\}$, 
$2 \leq m \leq M$.
By definition $\I^{(m)} \cap \I^{(m')}= \emptyset, 1 \leq m \not= m' \leq M$.
Put $n=\sum_{m=1}^{M} N_m$ and
$
\tau_j=\sum_{m=1}^{M} t_m \1(j \in \I^{(m)}), 
1 \leq j \leq n.
$
That is, 
$$
\tau_j=t_m \Longleftrightarrow j \in \I^{(m)}, \quad
1 \leq j \leq n, \quad 1 \leq m \leq M.
$$
Then the integrand in the LHS of (\ref{eqn:E=det})
is simply written as
\begin{equation}
\prod_{j=1}^n \chi_{\tau_j}(x_j)
\det_{1 \leq j, k \leq n} [ 
\mbK_{\xi}(\tau_j, x_j; \tau_k, x_k)],
\label{eqn:det_Z1}
\end{equation}
and the integral 
$\int_{\prod_{m=1}^{M} \W_{N_m}(S)} 
\prod_{m=1}^M d \x_{N_m}^{(m)} (\cdot)$
can be replaced by
$\{\prod_{m=1}^{M} N_m ! \}^{-1}$
$\int_{S^n} d \x \, (\cdot)$.
By the formula (\ref{eqn:det2}) of determinant,
(\ref{eqn:det_Z1}) is expressed as
\begin{eqnarray}
&& \sum_{\sigma \in \cS_n} (-1)^{n-\ell(\sigma)}
\prod_{\lambda=1}^{\ell(\sigma)}
\prod_{j=1}^{q_{\lambda}}
\chi_{\tau_{c_{\lambda}(j)}}(x_{c_{\lambda}(j)})
\mbK_{\xi}(\tau_{c_{\lambda}(j)}, x_{c_{\lambda}(j)};
\tau_{c_{\lambda}(j+1)}, x_{c_{\lambda}(j+1)})
\nonumber\\
&=& 
\sum_{\sigma \in \cS_n} (-1)^{n-\ell(\sigma)}
\prod_{\lambda=1}^{\ell(\sigma)}
\prod_{j=1}^{q_{\lambda}}
\chi_{\tau_{c_{\lambda}(j)}}(x_{c_{\lambda}(j)})
\Big\{
\cG_{\xi}
(\tau_{c_{\lambda}(j)}, x_{c_{\lambda}(j)}; \tau_{c_{\lambda}(j+1)}, x_{c_{\lambda}(j+1)})
\nonumber\\
&& \qquad \qquad 
-\1(\tau_{c_{\lambda}(j)} > \tau_{c_{\lambda}(j+1)})
p(\tau_{c_{\lambda}(j)} - \tau_{c_{\lambda}(j+1)},
x_{c_{\lambda}(j)} | x_{c_{\lambda}(j+1)} )
\Big\},
\label{eqn:star2}
\end{eqnarray}
where with (\ref{eqn:K1}) we have set 
\begin{eqnarray}
\cG_{\xi}(s, x;  t, y)
&\equiv& \mbK_{\xi}(s,x;t,y)+\1(s>t) p(s-t,x|y)
\nonumber\\
&=& \int_{S} \xi(dv) p(s, x|v) \cM_{\xi}^{v}(t,y),
\quad (s,t) \in [0, \infty)^2, \quad (x,y) \in S^2.
\label{eqn:G1}
\end{eqnarray}

We will perform binomial expansions in (\ref{eqn:star2});
for each $(\lambda, j)$, we choose a term
$\cG_{\xi}$
or a term
$-\1(\tau_{c_{\lambda}(j)} > \tau_{c_{\lambda}(j+1)}) p$.
Here we have to note that it is not a usual binomial expansion
in the sense that, since the latter term involves
an indicator function $\1(\tau_{c_{\lambda}(j)} > \tau_{c_{\lambda}(j+1)})$,
even if the latter term is chosen in the expansion,
nothing remains when
$\tau_{c_{\lambda}(j)}  \leq \tau_{c_{\lambda}(j+1)}$.
It means that, if the condition $\tau_{c_{\lambda}(j)} > \tau_{c_{\lambda}(j+1)}$
is not satisfied, we have no other choice than taking the first term.
Moreover, for each cyclic permutation labeled by $\lambda$, 
the first term has to be chosen at least once
in the $q_{\lambda}$-times expansion.
Otherwise, it vanishes, since as a matter of course
the conditions $\tau_{c_{\lambda}(j)}> \tau_{c_{\lambda}(j+1)},
1 \leq \forall j \leq q_{\lambda}$ cannot be satisfied 
in a cyclic permutation $\sc_{\lambda}$. 
For each cyclic permutation $\sc_{\lambda}$,
let $\bM_{\lambda}$ be the set of first terms and 
$\bM_{\lambda}^{\rm c}=\{\sc_{\lambda}\} \setminus \bM_{\lambda}$ 
be the set of non-vanishing second terms
in this restricted binomial expansion.
(Note again that $\bM_{\lambda} \not= \emptyset$.) 
If we define a subset of $\{\sc_{\lambda}\}$,
\begin{equation}
\cC(\sc_{\lambda}) =
\Big\{ c_{\lambda}(j) \in \{\sc_{\lambda}\} : 
\tau_{c_{\lambda}(j)} > \tau_{c_{\lambda}(j+1)} \Big\}, 
\label{eqn:cC1}
\end{equation}
then 
$\{\sc_{\lambda}\} \setminus \cC(\sc_{\lambda})
\subset \bM_{\lambda}$ 
and $\bM_{\lambda}^{\rm c} \subset \cC(\sc_{\lambda})$.
It implies that, if we put
\begin{eqnarray}
G(\sc_{\lambda}, \bM_{\lambda}) &=& 
\int_{S^{\{\sc_{\lambda}\}}} \prod_{j=1}^{q_{\lambda}} \ \bigg\{
dx_{c_{\lambda}(j)}  \  \chi_{{\tau}_{c_{\lambda}(j)}}(x_{c_{\lambda}(j)})
p( \tau_{c_{\lambda}(j)} - \tau_{c_{\lambda}(j+1)},
x_{c_{\lambda}(j)} | x_{c_{\lambda}(j+1)} )
^{\1(c_{\lambda}(j) \in \bM_{\lambda}^{\rm c})}
\nonumber\\
&& \qquad \qquad  \times
\cG_{\xi}
({\tau}_{c_{\lambda}(j)}, x_{c_{\lambda}(j)}; {\tau}_{c_{\lambda}(j+1)}, x_{c_{\lambda}(j+1)})
^{\1(c_{\lambda}(j) \in \bM_{\lambda})}\bigg\},
\label{eqn:G}
\end{eqnarray}
then, as a result of this expansion, 
the LHS of (\ref{eqn:E=det}) is expressed as
\begin{equation}
\frac{1}{\prod_{m=1}^{M} N_m !} 
\sum_{\sigma \in \cS_n} (-1)^{n-\ell(\sigma)}
\prod_{\lambda=1}^{\ell(\sigma)}
\sum_{\substack{\bM_{\lambda}: \\
\{\sc_{\lambda}\} \setminus \cC(\sc_{\lambda}) 
\subset \bM_{\lambda} \subset \{\sc_{\lambda}\}}}
(-1)^{\sharp \bM^{\rm c}_{\lambda}}
G(\sc_{\lambda}, \bM_{\lambda}).
\label{eqn:LHS2}
\end{equation}

\vskip 0.3cm
\noindent
{\it Step 2: Application of Fubini's theorem}
\vskip 0.3cm
Put
\begin{eqnarray}
&& F(\{x_{c_{\lambda}(k)}:c_{\lambda}(k) \in \bM_{\lambda}^{\rm c} \}) 
\nonumber\\
&& \qquad
= \int_{S^{\bM_{\lambda}}}
\prod_{j: c_{\lambda}(j) \in \bM_{\lambda}}\Big\{ dx_{c_{\lambda}(j)} 
\ \chi_{{\tau}_{c_{\lambda}(j)}}(x_{c_{\lambda}(j)})
\cG_{\xi}
({\tau}_{c_{\lambda}(j)}, x_{c_{\lambda}(j)}; {\tau}_{c_{\lambda}(j+1)}, x_{c_{\lambda}(j+1)}) \Big\}
\nonumber\\
&& \qquad \qquad \times
\prod_{k: c_{\lambda}(k) \in \bM_{\lambda}^{\rm c}}
p( \tau_{c_{\lambda}(k)} - \tau_{c_{\lambda}(k+1)},
x_{c_{\lambda}(k)} | x_{c_{\lambda}(k+1)} ),
\label{eqn:F}
\end{eqnarray}
which is the integral only over $S^{\bM_{\lambda}}$.
Then (\ref{eqn:G}) is obtained by performing the integral of (\ref{eqn:F}) 
over $S^{\bM^{\rm c}_{\lambda}}=S^{\{\sc_{\lambda}\}}
\setminus S^{\bM_{\lambda}}$, 
\begin{equation}
G(\sc_{\lambda}, \bM_{\lambda}) =
\int_{S^{\bM_{\lambda}^{\rm c}}} 
\prod_{k:c_{\lambda}(k) \in \bM_{\lambda}^{\rm c}} \Big\{ dx_{c_{\lambda}(k)} 
\chi_{{\tau}_{c_{\lambda}(k)}}(x_{c_{\lambda}(k)}) \Big\}
F(\{x_{c_{\lambda}(k)}:c_{\lambda}(k) \in \bM_{\lambda}^{\rm c} \}).
\label{eqn:G2}
\end{equation}
In (\ref{eqn:F}), we use the integral expression
(\ref{eqn:G1}) for each
$\cG_{\xi}
({\tau}_{c_{\lambda}(j)}, x_{c_{\lambda}(j)}; {\tau}_{c_{\lambda}(j+1)}, x_{c_{\lambda}(j+1)})$
by letting the integral variable
be $v=v_{c_{\lambda}(j)}$, then
\begin{eqnarray}
&& 
F(\{x_{c_{\lambda}(k)}:c_{\lambda}(k) \in \bM_{\lambda}^{\rm c}\}) 
\nonumber\\
&=& \int_{S^{ \bM_{\lambda}}}
\prod_{j: c_{\lambda}(j) \in \bM_{\lambda}} \xi (dv_{c_{\lambda}(j)})
\int_{S^{ \bM_{\lambda}}}
\prod_{j: c_{\lambda}(j) \in \bM_{\lambda}} 
\Big\{dx_{c_{\lambda}(j)}
p(\tau_{c_{\lambda}(j)}, x_{c_{\lambda}(j)}|v_{c_{\lambda}(j)})
\chi_{{\tau}_{c_{\lambda}(j)}}(x_{c_{\lambda}(j)}) \Big\}
\nonumber\\
&\times&
\prod_{j:c_{\lambda}(j) \in \bM_{\lambda}} 
\cM_{\xi}^{v_{c_{\lambda}(j)}}({\tau}_{c_{\lambda}(j+1)}, x_{c_{\lambda}(j+1)})
\prod_{k: c_{\lambda}(k) \in \bM_{\lambda}^{\rm c}} 
p({\tau}_{c_{\lambda}(k)}-{\tau}_{c_{\lambda}(k+1)}, 
x_{c_{\lambda}(k)}| x_{c_{\lambda}(k+1)}).
\nonumber
\end{eqnarray}
Let $\v=(v_{c_{\lambda}(j)}: c_{\lambda}(j) \in \bM_{\lambda})$.
Since $p(t, x|v_{c_{\lambda}(j)}) dx$ gives the probability
$\rP_{v_{c_{\lambda}(j)}}[Y_{c_{\lambda}(j)}(t) \in dx]$, the above 
expresses
\begin{eqnarray}
&&
\int_{S^{ \bM_{\lambda}}}
\prod_{j: c_{\lambda}(j) \in \bM_{\lambda}} \xi (dv_{c_{\lambda}(j)})
\rE_{\v} \Bigg[ 
\prod_{j:c_{\lambda}(j) \in \bM_{\lambda}} 
\bigg\{ \chi_{{\tau}_{c_{\lambda}(j)}}(Y_{c_{\lambda}(j)}({\tau}_{c_{\lambda}(j)}))
\nonumber\\
&& \qquad \qquad \qquad \times
\cM_{\xi}^{v_{c_{\lambda}(j)}}({\tau}_{c_{\lambda}(j+1)}, Y_{c_{\lambda}(j+1)}
({\tau}_{c_{\lambda}(j+1)}))^{\1(c_{\lambda}(j+1) \in \bM_{\lambda})}
\nonumber\\
&& \qquad \qquad \qquad  \times 
\cM_{\xi}^{v_{c_{\lambda}(j)}}
({\tau}_{c_{\lambda}(j+1)}, x_{c_{\lambda}(j+1)})^{\1(c_{\lambda}(j+1) 
\in \bM_{\lambda}^{\rm c})}\bigg\}
\nonumber\\
&& \quad \times
\prod_{k: c_{\lambda}(k) \in \bM_{\lambda}^{\rm c}} 
\bigg\{ p({\tau}_{c_{\lambda}(k)}-{\tau}_{c_{\lambda}(k+1)}, 
x_{c_{\lambda}(k)} | Y_{c_{\lambda}(k+1)}({\tau}_{c_{\lambda}(k+1)})
)^{\1(c_{\lambda}(k+1) \in \bM_{\lambda})} \nonumber\\
&& \qquad \qquad \qquad \qquad \times 
p({\tau}_{c_{\lambda}(k)}-{\tau}_{c_{\lambda}(k+1)}, 
x_{c_{\lambda}(k)} | x_{c_{\lambda}(k+1)})
^{\1(c_{\lambda}(k+1) \in \bM_{\lambda}^{\rm c})} \bigg\} \Bigg].
\nonumber
\end{eqnarray}
Using Fubini's theorem, (\ref{eqn:G2}) is given by
\begin{eqnarray}
&& \int_{S^{\bM_{\lambda}}}
\prod_{j:c_{\lambda}(j) \in \bM_{\lambda}} \xi(dv_{c_{\lambda}(j)})
\rE_{\v} \Bigg[ 
\prod_{j:c_{\lambda}(j) \in \bM_{\lambda}} 
\chi_{{\tau}_{c_{\lambda}(j)}}(Y_{c_{\lambda}(j)}({\tau}_{c_{\lambda}(j)}))
\nonumber\\
&& \qquad \qquad \times 
\prod_{j: c_{\lambda}(j), c_{\lambda}(j+1) \in \bM_{\lambda}}
\cM_{\xi}^{v_{c_{\lambda}(j)}}({\tau}_{c_{\lambda}(j+1)}, 
Y_{c_{\lambda}(j+1)}({\tau}_{c_{\lambda}(j+1)}))
\nonumber\\
&& \quad \times 
\int_{S^{\bM_{\lambda}^{\rm c}}} 
\prod_{k: c_{\lambda}(k) \in \bM_{\lambda}^{\rm c}} 
\Big\{ dx_{c_{\lambda}(k)} \chi_{\tau_{c_{\lambda}(k)}}
(x_{c_{\lambda}(k)}) \Big\}
\nonumber\\
&& \qquad \qquad \times 
\prod_{k: c_{\lambda}(k) \in 
\bM_{\lambda}^{\rm c}, c_{\lambda}(k+1) \in \bM_{\lambda}}
p( \tau_{c_{\lambda}(k)} - \tau_{c_{\lambda}(k+1)},
x_{c_{\lambda}(k)} | Y_{c_{\lambda}(k+1)}(\tau_{c_{\lambda}(k+1)}))
\nonumber\\
&& \qquad \qquad \times
\prod_{k: c_{\lambda}(k), c_{\lambda}(k+1) \in \bM_{\lambda}^{\rm c}}
p({\tau}_{c_{\lambda}(k)}-{\tau}_{c_{\lambda}(k+1)}, 
x_{c_{\lambda}(k)} | x_{c_{\lambda}(k+1)})
\nonumber\\
&& \qquad \qquad \times
\prod_{j: c_{\lambda}(j) \in \bM_{\lambda}, 
c_{\lambda}(j+1) \in \bM_{\lambda}^{\rm c}}
\cM_{\xi}^{v_{c_{\lambda}(j)}}
(\tau_{c_{\lambda}(j+1)}, x_{c_{\lambda}(j+1)})
\Bigg].
\nonumber
\end{eqnarray}

We will perform integration over $x_{c_{\lambda}(k)}$'s
for $c_{\lambda}(k) \in \bM^{\rm c}_{\lambda}$
before taking the expectation $\rE_{\v}$.
For each $1 \leq j \leq q_{\lambda}$ such that
$c_{\lambda}(j) \in \bM_{\lambda}$, we define
\begin{eqnarray}
\oj &=& \min\{k > j : c_{\lambda}(k) \in \bM_{\lambda} \}
\nonumber\\
\uj &=& \max\{k < j : c_{\lambda}(k) \in \bM_{\lambda} \}. 
\label{eqn:ojuj}
\end{eqnarray} 
Then, for each $j$ such that $c_{\lambda}(j) \in \bM_{\lambda}$,
what we want to calculate is 
\begin{eqnarray}
&& \chi_{{\tau}_{c_{\lambda}(j)}}(Y_{c_{\lambda}(j)}({\tau}_{c_{\lambda}(j)}))
\Big\{ \prod_{k=\uj+1}^{j-1}
\int_{S} dx_{c_{\lambda}(k)} \chi_{\tau_{c_{\lambda}(k)}}
(x_{c_{\lambda}(k)}) 
\Big\}
\nonumber\\
&& \qquad \quad \times
p(\tau_{c_{\lambda}(j-1)}-\tau_{c_{\lambda}(j)}, 
x_{c_{\lambda}(j-1)} |Y_{c_{\lambda}(j)}(\tau_{c_{\lambda}(j)}))
\nonumber\\
&& \qquad \quad \times
\prod_{l=\uj+2}^{j-1} 
p(\tau_{c_{\lambda}(l-1)}-\tau_{c_{\lambda}(l)},
x_{c_{\lambda}(l-1)} | x_{c_{\lambda}(l)})
\cM_{\xi}^{v_{c_{\lambda}(\, \uj \,)}}
(\tau_{c_{\lambda}(\uj+1)}, x_{c_{\lambda}(\uj+1)}),
\nonumber
\end{eqnarray}
which is the multiple integral over $x_{c_{\lambda}(k)}$'s
with indices in interval $\uj+1 \leq k \leq j-1$.
Since it is written as
\begin{eqnarray}
&& \prod_{k=\uj+1}^{j-1} 
\int_{S} d x_{c_{\lambda}(k)} \, \cM_{\xi}^{v_{c_{\lambda}(\uj)}}
(\tau_{c_{\lambda}(\uj+1)}, x_{c_{\lambda}(\uj+1)})
\nonumber\\
&& \quad \times \chi_{\tau_{c_{\lambda}(\uj+1)}}(x_{c_{\lambda}(\uj+1)})
p(\tau_{c_{\lambda}(\uj+1)}-\tau_{c_{\lambda}(\uj+2)}, x_{c_{\lambda}(\uj+1)} |x_{c_{\lambda}(\uj+2)})
\nonumber\\
&& \quad \times \chi_{\tau_{c_{\lambda}(\uj+2)}}(x_{c_{\lambda}(\uj+2)})
p(\tau_{c_{\lambda}(\uj+2)}-\tau_{c_{\lambda}(\uj+3)}, x_{c_{\lambda}(\uj+2)} |x_{c_{\lambda}(\uj+3)})
\nonumber\\
&& \quad \times \cdots
\nonumber\\
&& \quad \times \chi_{\tau_{c_{\lambda}(j-2)}}(x_{c_{\lambda}(j-2)})
p(\tau_{c_{\lambda}(j-2)}-\tau_{c_{\lambda}(j-1)}, x_{c_{\lambda}(j-2)} |x_{c_{\lambda}(j-1)})
\nonumber\\
&& \quad \times \chi_{\tau_{c_{\lambda}(j-1)}}(x_{c_{\lambda}(j-1)})
p(\tau_{c_{\lambda}(j-1)}-\tau_{c_{\lambda}(j)}, x_{c_{\lambda}(j-1)} |Y_{c_{\lambda}(j)}(\tau_{c_{\lambda}(j)}))
\chi_{\tau_{c(j)}}(Y_{c_{\lambda}(j)}(\tau_{c_{\lambda}(j)})),
\nonumber
\end{eqnarray}
it represents the conditional expectation of
\begin{eqnarray}
&&\prod_{k=\uj+1}^{j} \chi_{\tau_{c_{\lambda}(k)}}
(Y_{c_{\lambda}(j)}(\tau_{c_{\lambda}(k)}))
\cM_{\xi}^{v_{c_{\lambda}(\, \uj \,)}}
(\tau_{c_{\lambda}(\uj+1)}, Y_{c_{\lambda}(j)}(\tau_{c_{\lambda}(\uj+1)}))
\nonumber
\end{eqnarray}
with respect to $\rE_{\v}[\, \cdot \,|Y_{c_\lambda(j)}(\tau_{c_{\lambda}(j)})]$.
Then 
(\ref{eqn:G}) is equal to 
\begin{eqnarray}
&& G(\sc_{\lambda}, \bM_{\lambda})=\int_{S^{\bM_{\lambda}}} 
\prod_{j: c_{\lambda}(j) \in \bM_{\lambda}} \xi(dv_{c_{\lambda}(j)})
\nonumber\\
&& \times
\rE_{\v} \left[ 
\prod_{j: c_{\lambda}(j) \in \bM_{\lambda}}
\left\{\prod_{k=\uj+1}^{j}
\chi_{\tau_{c_{\lambda}(k)}}(Y_{c_{\lambda}(j)}(\tau_{c_{\lambda}(k)})) 
\cM_{\xi}^{v_{c_{\lambda}(\uj)}}
(\tau_{c_{\lambda}(\uj+1)}, Y_{c_{\lambda}(j)}(\tau_{c_{\lambda}(\uj+1)}))
\right\} \right]. 
\nonumber\\
\label{eqn:Gb}
\end{eqnarray}

\vskip 0.3cm
\noindent
{\it Step 3: Constructions of sub-cyclic-permutation $\widehat{\sigma}$
and subsets $\{\bJ^{\lambda}_{m}\}$ }
\vskip 0.3cm
Using only the entries of $\bM_{\lambda} \not= \emptyset$,
we uniquely define a cyclic permutation $\widehat{\sc_{\lambda}}$ 
as a subsystem of $\sc_{\lambda}$ as follows. 
We write $\widehat{q_{\lambda}}=\sharp\{\widehat{\sc_{\lambda}}\}=\sharp \bM_{\lambda} \geq 1$.
Let $j_1=\min \{1 \leq j \leq q_{\lambda}: c_{\lambda}(j) 
\in \bM_{\lambda}\}$ and define $\widehat{c_{\lambda}}(1)=c_{\lambda}(j_1)$.
If $\widehat{q_{\lambda}} \geq 2$, we put
$j_{k+1}=\overline{j_k}, 1 \leq k \leq \widehat{q_{\lambda}}-1$.
Then $\widehat{\sc_{\lambda}} 
=(\widehat{c_{\lambda}}(1) \widehat{c_{\lambda}}(2) 
\cdots \widehat{c_{\lambda}}(\widehat{q_{\lambda}}))$
is defined by $(c_{\lambda}(j_1) c_{\lambda}(j_2) 
\cdots c(j_{\widehat{q_{\lambda}}}))$.

Moreover, 
we decompose the set $\bM_{\lambda}$ into $M$ subsets,
$\bM_{\lambda}=\bigcup_{m=1}^{M} \bJ_m^{\lambda}$,
by letting 
\begin{equation}
\bJ_m^{\lambda}=\bJ_m^{\lambda}(\sc_{\lambda}, \bM_{\lambda})=
\Big\{ c_{\lambda}(j) \in \bM_{\lambda} :
\uj < ^{\exists}k \leq j, \,
\mbox{s.t.} \, \tau_{c_{\lambda}(k)}=t_m \Big\}, \quad 1 \leq m \leq M.
\label{eqn:bJ1}
\end{equation}
By this definition,
\begin{equation}
\prod_{j: c_{\lambda}(j) \in \bM_{\lambda}}
\prod_{k=\uj+1}^{j}
\chi_{\tau_{c_{\lambda}(k)}}(Y_{c_{\lambda}(j)}(\tau_{c_{\lambda}(k)})) 
= \prod_{m=1}^M \prod_{j_m \in \bJ^{\lambda}_m}
\chi_{t_m}(Y_{j_m}(t_m)).
\label{eqn:eqchi}
\end{equation}
The subsets $\{\bJ_{m}^{\lambda}\}$ have the following properties.
Assume that $c_{\lambda}(j) \in \bM_{\lambda}$.  
\begin{description}
\item{(i)} \quad
By definition (\ref{eqn:bJ1}), 
if $\tau_{c_{\lambda}(j)}=t_m$, 
$c_{\lambda}(j) \in \bJ_m^{\lambda}$.
That is, $\I^{(m)} \cap \bM_{\lambda} \subset \bJ_m^{\lambda}, 1 \leq m \leq M$.

\item{(ii)} \quad
Assume $\tau_{c_{\lambda}(j)}=t_m$ and 
$\uj < j-2$. Then, by definition of $\uj$ given by (\ref{eqn:ojuj}),
for all $k$ such that
$\uj+1 \leq k \leq j-1$, 
$c_{\lambda}(k) \in \bM_{\lambda}^{\rm c} \subset \cC(\sc_{\lambda})$
and thus 
\begin{equation}
\tau_{c_{\lambda}(k)} > \tau_{c_{\lambda}(k+1)} \geq \tau_{c_{\lambda}(j)}=t_m.
\label{eqn:tau_ineq}
\end{equation}
Therefore, $c_{\lambda}(j)$ can belong to plural
$\bJ_{m'}^{\lambda}$'s with $m' \geq m$.
That is, $\bJ_m^{\lambda} \cap \bJ_{m'}^{\lambda} \not= \emptyset,
m \not= m'$, in general.

\item{(iii)} \quad
If $\tau_{c_{\lambda}(j)}=t_m$ with $m \geq 2$, 
then $c_{\lambda}(j) \notin \bJ_1^{\lambda}$ by (\ref{eqn:tau_ineq}).
It implies that, if $c_{\lambda}(j) \in \bJ_1^{\lambda}$,
then $\tau_{c_{\lambda}(j)}=t_1$. That is,
$$
\bJ_1^{\lambda}=\I^{(1)} \cap \bM_{\lambda}
= \I^{(1)} \cap \{\sc_{\lambda}\}.
$$
Here the second equality comes from the definition 
of (\ref{eqn:cC1});
if $\tau_{c_{\lambda}(j)}=1$, since it is the minimal value of 
$\tau_{c_{\lambda}(j)}$, 
$c_{\lambda}(j) \notin \cC(\sc_{\lambda})$.

\item{(iv)} \quad
If $\tau_{c_{\lambda}(j)}=t_m$ with $m \geq 3$, then 
$c_{\lambda}(j) \notin \bJ_{m'}^{\lambda}, 1 \leq m' \leq m-1$
by (\ref{eqn:tau_ineq}).
It implies that, if 
$c_{\lambda}(j) \in \bJ_{m}^{\lambda}, 2 \leq m \leq M$,
then $\tau_{c_{\lambda}(j)} \leq t_m$. That is,
$$
\bJ_{m}^{\lambda} \subset \bigcup_{k=1}^{m} \I^{(k)},
\quad 2 \leq m \leq M.
$$
Moreover, if $c_{\lambda}(j) \in \bJ_m^{\lambda}$
and $\tau_{c_{\lambda}(j)} = t_k$ with $k < m$,
by definition (\ref{eqn:bJ1}),
$c_{\lambda}(j) \in \bJ_k^{\lambda}$. That is,
$$
\bJ_{m}^{\lambda} \cap \I^{(k)} \subset \bJ_{k}^{\lambda}, \quad
1 \leq k < m \leq M.
$$
\end{description}

By definition of $\widehat{\sc}_{\lambda}$, 
$\prod_{j: c_{\lambda}(j) \in \bM_{\lambda}} \xi(d v_{c_{\lambda}(j)})$
in (\ref{eqn:Gb}) is written as
$\prod_{j=1}^{\widehat{q_{\lambda}}} \xi(d v_{\widehat{c_{\lambda}}(j)})$
and $\v$ is given by
$\v=(v_{\widehat{c_{\lambda}}(1)}, \dots, v_{\widehat{c_{\lambda}}( \widehat{q_{\lambda}})})$.
By (\ref{eqn:tau_ineq}), $\tau_{c_{\lambda}(\uj+1)} \geq \tau_{c_{\lambda}(k)},
\uj+1 \leq \forall k \leq j$.
By the setting, $\tau_{c_{\lambda}(\uj+1)} \leq t_M \leq T$.
Then, since $\cM_{\xi}^v$ is a continuous-time martingale, 
$\cM_{\xi}^{v_{c_{\lambda}(\uj)}}
(\tau_{c_{\lambda}(\uj+1)}, Y_{c_{\lambda}(j)}(\tau_{c_{\lambda}(\uj+1)}))$
can be replaced by
$\cM_{\xi}^{v_{c_{\lambda}(\uj)}}(T, Y_{c_{\lambda}(j)}(T))$,
for each $1 \leq j \leq \widehat{q_{\lambda}}$ in (\ref{eqn:Gb}).
Again, by definition of $\widehat{\sc_{\lambda}}$,
$$
\prod_{j: c_{\lambda}(j) \in \bM_{\lambda}}
\cM_{\xi}^{v_{c_{\lambda}(\uj)}}(T, Y_{c_{\lambda}(j)}(T))
= \prod_{j=1}^{\widehat{q_{\lambda}}}
\cM_{\xi}^{v_{\widehat{c_{\lambda}}(j)}}
(T, Y_{\widehat{c_{\lambda}}(j+1)}(T)).
$$
Therefore, (\ref{eqn:Gb}) is rewritten as
\begin{equation}
G(\sc_{\lambda}, \bM_{\lambda})=
\int_{S^{\bM_{\lambda}}} \prod_{j=1}^{\widehat{q_{\lambda}}}
\xi(dv_{\widehat{c_{\lambda}}(j)}) \
\rE_{\v} \left[ 
\prod_{m=1}^{M} \prod_{j_m \in \bJ_m^{\lambda}}
\chi_{t_m}(Y_{j_m}(t_m))
\prod_{j=1}^{\widehat{q_{\lambda}}}
\cM_{\xi}^{v_{\widehat{c_{\lambda}}(j)}} 
(T, Y_{\widehat{c_{\lambda}}(j+1)}(T)) \right],
\label{eqn:G4}
\end{equation}
where
$\v=(v_{\widehat{c_{\lambda}}(1)}, \dots, v_{\widehat{c_{\lambda}}(\widehat{q_{\lambda}})})$.

Let $\bM \equiv \bigcup_{\lambda=1}^{\ell(\sigma)}
\bM_{\lambda}$. 
Since
$n-\sum_{\lambda=1}^{\ell(\sigma)} \sharp \bM_{\lambda}^{\rm c}
=\sharp \bM$,
the LHS of (\ref{eqn:E=det}), which is expanded as 
(\ref{eqn:LHS2}) with (\ref{eqn:G4}), becomes now
\begin{eqnarray}
&& \frac{1}{\prod_{m=1}^{M} N_m !} 
\sum_{\sigma \in \cS_n} 
\sum_{\substack{\bM : \\
\I_n \setminus \bigcup_{\lambda=1}^{\ell(\sigma)}
\cC(\sc_{\lambda}) \subset \bM \subset \I_n}}
(-1)^{\sharp \bM -\ell(\sigma)}
\int_{S^{\bM}}
\prod_{\lambda=1}^{\ell(\sigma)}
\prod_{j=1}^{\widehat{q_{\lambda}}}
\xi(dv_{\widehat{c_{\lambda}}(j)})
\nonumber\\
&& \qquad \times 
\rE_{\v} \left[ \prod_{\lambda=1}^{\ell(\sigma)}
\left\{ \prod_{m=1}^{M}
\prod_{j_m \in \bJ_m^{\lambda}} 
\chi_{t_m}(Y_{j_m}(t_m)) 
\prod_{j=1}^{\widehat{q_{\lambda}}}
\cM_{\xi}^{v_{\widehat{c_{\lambda}}(j)}}
(T, Y_{\widehat{c_{\lambda}}(j+1)}(T)) 
\right\} \right].
\label{eqn:E=det2}
\end{eqnarray}
We define 
\begin{equation}
\widehat{\sigma} \equiv 
\widehat{\sc_1} \widehat{\sc_2} \cdots
\widehat{\sc_{\ell(\sigma)}}
\label{eqn:sigma_hat_prod}
\end{equation}
and  
$$
\bJ_m \equiv \bigcup_{\lambda=1}^{\ell(\sigma)} \bJ_m^{\lambda},
\quad 1 \leq m \leq M.
$$
Note that $\ell(\widehat{\sigma})=\ell(\sigma)$,
since $\bM_{\lambda}\not= \emptyset$ and thus
$\widehat{\sc_{\lambda}}$ is not empty for all $1 \leq \lambda \leq \ell(\sigma)$.
Due to the properties of $\bJ^{\lambda}_m, 1 \leq \lambda \leq \ell(\sigma)$
given below (\ref{eqn:eqchi}), 
the obtained $(\bJ_m)_{m=1}^{M}$'s 
form a collection of series of index sets
satisfying the following conditions,
which we write as ${\cal J}(\{N_m\}_{m=1}^{M})$:
\begin{eqnarray}
&&
\bJ_1=\I^{(1)},  \quad
\bJ_m \subset \bigcup_{k=1}^m \I^{(k)}
\quad \mbox{for} \quad 2 \leq m \leq M,
\nonumber\\
&& 
\bJ_m \cap \I^{(k)} \subset \bJ_{k} 
\quad \mbox{for} \quad
1 \leq k < m \leq M,  \quad \mbox{and}
\nonumber\\
&& 
\sharp \bJ_m=N_m \quad \mbox{for} \quad 1 \leq m \leq M.
\nonumber
\end{eqnarray}
By definition
\begin{equation}
\bM=\bigcup_{m=1}^M \bJ_m.
\label{eqn:Mtotal}
\end{equation}

\vskip 0.3cm
\noindent
{\it Step 4: Combinatorics arguments}
\vskip 0.3cm
In the previous steps, given an index set with $M$ levels
({\it i.e.} $M$ different times $t_m, 1 \leq m \leq M$),
$\I_n=\bigcup_{m=1}^M \I^{(m)}$
and a permutation $\sigma \in \cS_n$,
we gave a procedure to obtain a sub-index-set $\bM$,
a sub-permutation $\widehat{\sigma} \in \cS(\bM)$,
and a series of index sets $(\bJ_m)_{m=1}^M \in {\cal J}(\{N_m\}_{m=1}^{M})$.
For obtained $(\bJ_m)_{m=1}^M$, let
$$
\bA_1=\emptyset \quad \mbox{and} \quad
\bA_m=\bJ_m \cap \bigcup_{k=1}^{m-1} \I^{(k)}
=\bJ_m \cap \bigcup_{k=1}^{m-1} \bJ_k, \quad 
2 \leq m \leq M.
$$
We set $A_m=\sharp \bA_m, 1 \leq m \leq M$.
In the reduction from $\sigma \in \cS_n$ to 
$\widehat{\sigma} \in \cS(\bM)$ and $(\bJ_m)_{m=1}^M$
with (\ref{eqn:Mtotal}), the indices in the set $\bA_m$ are eliminated
from the original index set $\I^{(m)}$ at each level $1 \leq m \leq M$.
Then 
$$
\sharp \bM=\sum_{m=1}^M (N_m-A_m), 
$$
since $\sharp \I^{(m)}=N_m, 1 \leq m \leq M$.
Now we consider reversing the procedure.
Given $\widehat{\sigma} \in \cS(\bM)$ and $(\bJ_m)_{m=1}^M$
with (\ref{eqn:Mtotal}), if we represent $\widehat{\sigma}$ as
(\ref{eqn:sigma_hat_prod}) and insert the indices in $\bA_m$ at each level
$1 \leq m \leq M$ into cycles 
$\widehat{\sc_{\lambda}}, 1 \leq \lambda \leq \widehat{q_{\lambda}}$
by perfectly tracing back the procedure,
we will obtain the original permutation $\sigma \in \cS_n$.
But if we exchange the indices within each set $\bA_m, 1 \leq m \leq M$,
when we insert them into $\widehat{\sc_{\lambda}}$'s, 
we will obtain a different permutation, say
$\sigma' \in \cS_n$ (keeping $\ell(\sigma')=\ell(\sigma)$).
It implies that 
the number of distinct $\sigma$'s in $\cS_n$ which give
the same $\widehat{\sigma}$ and
$(\bJ_m)_{m=1}^{M}$ by the present reduction is given by
$\prod_{m=1}^{M} A_m !$, where $0! \equiv 1$.
Then (\ref{eqn:E=det2}) is equal to 
\begin{eqnarray}
&& 
\frac{1}{\prod_{m=1}^{M} N_m !}
\sum_{\substack{\bM: \\
\max_m\{N_m\} \leq \sharp \bM \leq N}}
\sum_{\substack{(\bJ_m)_{m=1}^{M} \subset
{\cal J}(\{N_m\}_{m=1}^{M}): \\
\bigcup_{m=1}^{M} \bJ_m=\bM
}}
\prod_{m=1}^{M} A_m !
\sum_{\widehat{\sigma} \in \cS(\bM)}
(-1)^{\sharp \bM-\ell(\widehat{\sigma})}
\nonumber\\
&& \qquad \times
\sharp \bM !
\int_{\W_{\sharp \bM}(S)} \xi^{\otimes \bM} (d \v)
\rE_{\v} \left[
\prod_{m=1}^{M} \prod_{j_m \in \bJ_m}
\chi_{t_m}(Y_{j_m}(t_m))
\prod_{\lambda=1}^{\ell(\widehat{\sigma})}
\prod_{j=1}^{\widehat{q_{\lambda}}}
\cM_{\xi}^{v_{\widehat{ c_{\lambda} }(j)}}
(T, Y_{\widehat{c_{\lambda}}(j+1)}(T)) \right]
\nonumber\\
&&=
\sum_{\substack{\bM: \\
\max_m\{N_m\} \leq \sharp \bM \leq N}}
\sum_{\substack{(\bJ_m)_{m=1}^{M} \subset
{\cal J}(\{N_m\}_{m=1}^{M}): \\
\bigcup_{m=1}^{M} \bJ_m=\bM
}}
\sharp \bM ! \prod_{m=1}^{M} \frac{A_m!}{N_m!}
\nonumber\\
&& \qquad \times
\int_{\W_{\sharp \bM}(S)} \xi^{\otimes \bM} (d \v)
\rE_{\v} \left[
\prod_{m=1}^{M} \prod_{j_m \in \bJ_m}
\chi_{t_m}(Y_{j_m}(t_m))
\det_{j,k \in \bM}
\left[ \cM_{\xi}^{v_k} (T, Y_{j}(T)) \right] \right].
\label{eqn:AAA}
\end{eqnarray}

Let $1 \leq p \leq N$, 
$0 \leq A_m \leq N_m, 2 \leq m \leq M$, and $A_1 =0$.
Assume that
\begin{equation}
\sum_{m=1}^M (N_m-A_m)=p.
\label{eqn:restriction1}
\end{equation}
Consider
\begin{eqnarray}
&& \Lambda_1 =\left\{(\bJ_m)_{m=1}^{M} \subset
{\cal J}(\{N_{m}\}_{m=1}^{M}): 
\sharp \left(\bigcup_{m=1}^{M} \bJ_m \right) =p,
\right. 
\nonumber\\
&& \hskip 6cm \left.
\sharp \left(\bJ_m \cap \bigcup_{k=1}^{m-1} \bJ_k \right)
=A_m, 2 \leq m \leq M \right\},
\nonumber\\
&& \Lambda_2 = \left\{(\J_m)_{m=1}^{M} :
\sharp \J_m=N_m, 1 \leq m \leq M,
\bigcup_{m=1}^{M} \J_m= \I_{p},
\right. 
\nonumber\\
&& \hskip 6cm \left.
\sharp \left(\J_m \cap \bigcup_{k=1}^{m-1} \J_k \right)
=A_m, 2 \leq m \leq M \right\}.
\nonumber
\end{eqnarray}
Since $Y_j(\cdot)$'s are equivalent in probability law in $\rP_{\v}$, 
the integral in (\ref{eqn:AAA}) has the same value
for all $(\bJ_m)_{m=1}^{M} \in \Lambda_1$
with (\ref{eqn:Mtotal})
and it is also equal to
$$
\int_{\W_{p}(S)} \xi^{\otimes p} (d \v)
\rE_{\v} \left[
\prod_{m=1}^{M} \prod_{j_m \in \J_m}
\chi_{t_m}(Y_{j_m}(t_m))
\det_{j,k \in \I_{p}}
\left[ \cM_{\xi}^{v_k} (T, Y_j(T)) \right] \right]
$$
for all $(\J_{m})_{m=1}^{M} \in \Lambda_2$.

In $\Lambda_1$, $\bJ_1 \equiv \I^{(1)}$,
and for each set $\bJ_m$, $2 \leq m \leq M$,
$A_m$ elements in $\bJ_m$ are chosen from
$\bigcup_{k=1}^{m-1} \bJ_k$,
where $\sharp (\bigcup_{k=1}^{m-1}
\bJ_{k})=\sum_{k=1}^{m-1}(N_k-A_k)$,
and the remaining $N_m-A_m$ elements in $\bJ_m$
are from $\I^{(m)}$, where $\sharp \I^{(m)}=N_m$.
Then
$$
\sharp \Lambda_1=\prod_{m=2}^{M} 
\left\{ {\sum_{k=1}^{m-1} (N_k-A_k) \choose A_m}
{N_m \choose N_m-A_m} \right\}.
$$
In $\Lambda_2$, on the other hand,
$N_1$ elements in $\J_1$ is chosen from $\I_{p}$,
and then for each $2 \leq m \leq M$,
$A_m$ elements in $\J_m$ are chosen from
$\bigcup_{k=1}^{m-1} \J_k$, where
$\sharp (\bigcup_{k=1}^{m-1} \J_k)=
\sum_{k=1}^{m-1}(N_k-A_k)$, 
and the remaining $N_m-A_m$ elements in $\J_m$
are from $\I_{p} \setminus \bigcup_{k=1}^{m-1} \J_k$, 
where $\sharp(\I_{p} \setminus \bigcup_{k=1}^{m-1} \J_k)
=p-\sum_{k=1}^{m-1}(N_k-A_k)$.
Then
$$
\sharp \Lambda_2
={p \choose N_1} 
\prod_{m=2}^{M} \left\{
{\sum_{k=1}^{m-1} (N_k-A_k) \choose A_m}
{p-\sum_{k=1}^{m-1} (N_k-A_k) \choose N_m-A_m} \right\}.
$$
By (\ref{eqn:restriction1}),
we see 
$\sharp \Lambda_2/\sharp \Lambda_1=
p ! \prod_{m=1}^{M} A_m !/N_m!$.
Then (\ref{eqn:AAA}) is equal to the RHS of (\ref{eqn:E=det})
and the proof is completed.
\qed

\subsection{Proof of Theorem \ref{thm:DM_det} \label{sec:proof1_Theorem}}
By (\ref{eqn:Xi1}),
the moment generating function (\ref{eqn:GF1}) is 
written in the form
$$
\Psi_{\xi}^{\t}[\f]
=\E_{\xi} \left[ \prod_{m=1}^N \prod_{j=1}^N
\{1+\chi_{t_m}(X_j(t_m)) \} \right].
$$
By assumption of the theorem, it has DMR associated with $(Y,\cM_{\xi})$,
$$
\Psi_{\xi}^{\t}[\f]
=\rE_{\u} \left[
\prod_{m=1}^N \prod_{j=1}^N
\{1+\chi_{t_m}(Y_j(t_m)) \} 
\cD_{\xi}(T, \Y(T)) \right].
$$
Then Lemma \ref{thm:Fredholm} gives
a Fredholm determinant expression
to this as
$$
\Psi_{\xi}^{\t}[\f]
=\mathop{{\rm Det}}_
{\substack{
(s,t)\in \{t_1,t_2,\dots, t_M\}^2, \\
(x,y)\in S^2}
}
 \Big[\delta_{st} \delta_x(\{y\})
+ \mbK_{\xi}(s,x;t,y) \chi_{t}(y) \Big]
$$
with (\ref{eqn:K1}). 
By Definition \ref{thm:determinantal},
the proof is completed. \qed

\SSC{Preliminaries \label{sec:martingales}}

\subsection{Fundamental martingale-polynomials}
For $n \in \N_0$, we consider the monic
polynomials of degrees $n$ with time-dependent 
coefficients,
$$
m_n(t, x)=x^n + \sum_{j=0}^{n-1} c_n^{(j)}(t) x^j, \quad t \geq 0
$$
satisfying the conditions such that
$
m_n(0,x)=x^n,
$
and that, if we replace $x$ by the process
$Y(t), t \in [0, \infty)$, then they provide martingales.
If there exist such polynomials $\{m_n(t,x)\}_{n \in \N_0}$,
$t \in [0, \infty)$, $x \in S$,
we call them {\it fundamental martingale polynomials} (FMPs)
associated with the process $Y$.

\vskip 0.5cm
\noindent{\it BM and Hermite polynomials}
\vskip 0.3cm
Let $Y(t), t \in [0, \infty)$ be BM on $S=\R$.
The TD for BM on $\R$ is given by (\ref{eqn:p}),
$p=p_{\rm BM}$. 
For $n \in \N_0$, 
the Hermite polynomials of degrees $n \in \N_0$
are defined by
$$
H_n(x)=\sum_{j=0}^{[n/2]} (-1)^j 
\frac{n!}{j! (n-2j)!} (2x)^{n-2j}, 
\quad n \in \N_0, \quad x \in \R.
$$
The following is well-known 
(see, for instance, \cite{Sch00}).
\begin{lem}
\label{thm:Hermite1}
The FMPs
associated with BM on $\R$
are given by
$$
m_n(t, x)
=\left( \frac{t}{2} \right)^{n/2} 
H_n \left( \frac{x}{\sqrt{2t}} \right),
\quad t \geq 0, \quad n \in \N_0.
$$
\end{lem}
\noindent{\it Proof.} \,
Let
\begin{equation}
G_{\alpha}(t, y)
= e^{\alpha y-\alpha^2 t/2}, \quad \alpha \in \C.
\label{eqn:G_alpha}
\end{equation}
By (\ref{eqn:def_p}), we can show that, for $0 \leq s \leq t$, 
\begin{eqnarray}
\rE[ G_{\alpha}(t, Y(t)) | \cF(s)]
&=& \int_{\R} dy \, e^{\alpha y-\alpha^2 t/2} p_{\rm BM}(t-s, y|Y(s))
\nonumber\\
&=& \frac{e^{-Y(s)^2/2(t-s)-\alpha^2 t/2}}{\sqrt{2 \pi (t-s)}}
\int_{\R} dy \, e^{-(y^2-2Y(s)y)/2(t-s)+\alpha y}
\nonumber\\
&=& e^{\alpha Y(s)-\alpha^2 s/2} = G_{\alpha}(s, Y(s))
\quad \mbox{a.s.}
\label{eqn:G_martingale}
\end{eqnarray}
That is, $G_{\alpha}(t, Y(t)), t \in [0, \infty)$ is 
a continuous-time martingale.
The function (\ref{eqn:G_alpha}) is expanded using the Hermite polynomials as
\begin{equation}
G_{\alpha}(t,x)
=\sum_{n=0}^{\infty} \left( \frac{t}{2} \right)^{n/2}
H_n \left( \frac{x}{\sqrt{2t}} \right) 
\frac{\alpha^n}{n!}.
\label{eqn:GBM2}
\end{equation}
Then Lemma \ref{thm:Hermite1} is immediately concluded,
since we can easily confirm that $m_n(t,x), n \in \N_0$ 
are monic. 
\qed
\vskip 0.3cm

\vskip 0.5cm
\noindent{\it BESQ$^{(\nu)}$ and Laguerre polynomials}
\vskip 0.3cm
Next let $Y(t), t \in [0, \infty)$ be
BESQ$^{(\nu)}$ with $\nu>-1$ on $S=\R_+$.
For $n \in \N_0$, 
the Laguerre polynomials of degrees $n \in \N_0$ with 
index $\nu > -1$ 
are defined by
$$
L^{(\nu)}_n(x)=\sum_{j=0}^{n}(-1)^j 
\frac{\Gamma(n+\nu+1)}{\Gamma(\nu+j+1) (n-j)! j!} x^j,
\quad n \in \N_0, \quad x \in \R_+.
$$
\begin{lem}
\label{thm:Laguerre1}
Let $\nu > -1$. 
The FMPs
associated with BESQ$^{(\nu)}$ on $\R_+$
are given by
$$
m_n(t, x)
=(-1)^n n! (2t)^{n} L^{(\nu)}_n
\left( \frac{x}{2t} \right),
\quad t \geq 0, \quad n \in \N_0.
$$
\end{lem}
\noindent{\it Proof.} \,
For $\alpha \geq 0, t \geq 0, x \in \R_+, \nu > -1$, let
\begin{equation}
G_{\alpha}^{(\nu)}(t,x)
=\frac{e^{\alpha x/(1+2 \alpha t)}}{(1+2 \alpha t)^{\nu+1}}.
\label{eqn:Gnu1}
\end{equation}
For $0 \leq s \leq t$, we will calculate
$$
\rE[G_{\alpha}^{(\nu)}(t, Y(t))| \cF(s)]
=\frac{1}{(1+2 \alpha t)^{\nu+1}}
\int_{\R_+} dy \, e^{\alpha y/(1+2 \alpha t)} p^{(\nu)}(t-s, y|Y(s)) .
$$
If $Y(s)=0$,
\begin{eqnarray}
&& \int_{\R_+} dy \, e^{\alpha y/(1+2 \alpha t)} p^{(\nu)}(t-s, y|0) 
\nonumber\\
&& \qquad = \frac{1}{(2(t-s))^{\nu+1} \Gamma(\nu+1)}
\int_{\R_+} dy \, y^{\nu}
e^{-(1+2 \alpha s)y/\{2(t-s)(1+2 \alpha t)\}}
\nonumber\\
&& \qquad = \left( \frac{1+2 \alpha t}{1+2 \alpha s} \right)^{\nu+1}.
\nonumber
\end{eqnarray}
If $Y(s) >0$, 
\begin{eqnarray}
&& \int_{\R_+} dy \, e^{\alpha y/(1+2t \alpha)} p^{(\nu)}(t-s, y|Y(s)) 
\nonumber\\
&& \qquad = \frac{e^{-Y(s)/2(t-s)}}{2(t-s) Y(s)^{\nu/2}}
\int_{\R_+} dy \, y^{\nu}
e^{-(1+2 \alpha s)y/\{2(t-s)(1+2 \alpha t)\}}
I_{\nu} \left( \frac{\sqrt{Y(s) y}}{t-s} \right)
\nonumber\\
&& \qquad = \left( \frac{1+2 \alpha t}{1+2 \alpha s} \right)^{\nu+1}
e^{\alpha Y(s)/(1+2 \alpha s)} \quad \mbox{a.s.}
\nonumber
\end{eqnarray}
Here we used the integral formula
$$
\int_{\R_+} du \, u^{\nu/2} e^{-p^2 u} I_{\nu}(q \sqrt{u})
= \frac{q^{\nu}}{2^{\nu} p^{2(\nu+1)}} e^{q^2/4p},
\quad \nu> -1, \quad p \in \R, \quad q > 0,
$$
which is obtained from Eq.(4) of Section 13.3 in \cite{Wat44},
\begin{equation}
\int_{\R_+} dt \, t^{\nu+1} e^{-p^2 t^2} J_{\nu}(a t)
= \frac{a^{\nu}}{(2p^2)^{\nu+1}} e^{-a^2/4p},
\quad \Re \nu> -1, \quad |{\rm arg} \, p| < \frac{\pi}{4}, \quad a \in \C,
\label{eqn:int_form2}
\end{equation}
by setting $t=\sqrt{u}, a=iq$ and using the equality
\begin{equation}
I_{\nu}(z)=e^{-\nu \pi i /2} J_{\nu}(i z), \quad - \pi < {\rm arg} \, z < \frac{\pi}{2}.
\label{eqn:I_J}
\end{equation} 
Therefore, $G_{\alpha}^{(\nu)}(t, Y(t)), t \in [0, \infty)$
is a continuous-time martingale; 
$$
\rE[G_{\alpha}^{(\nu)}(t, Y(t)) | \cF(s)]
=G_{\alpha}^{(\nu)}(s, Y(s)) \quad \mbox{a.s. for $0 \leq s \leq t$}.
$$
The function (\ref{eqn:Gnu1}) is expanded by using Laguerre polynomials
as
$$
G_{\alpha}^{(\nu)}(t,x)
=\sum_{n=0}^{\infty} (-1)^n n! (2t)^n
L^{(\nu)}_n \left( \frac{x}{2 t} \right)
\frac{\alpha^n}{n!}.
$$
Then the lemma is proved.
\qed
\vskip 0.3cm

\subsection{Integral representations of martingales}

For each set of FMPs
$\{m_n(t,x): n \in \N_0\}$ associated with
the process $Y(t)$, here we want to
determine the integral transformation
of a function $f$
of the form
\begin{equation}
\cI[f(W) |(t, x)]
=\int_{S} dw \, f(cw) q(t, w|x)
\label{eqn:int1}
\end{equation}
with some constant $c \in \C$, 
such that it
satisfies the equalities
\begin{equation}
m_n(t, x)=\cI \left[ \left.
W^n \right|(t, x) \right],
\quad \forall n \in \N_0, \quad \forall t \in [0, \infty).
\label{eqn:int2}
\end{equation}
If so, 
given any polynomial $f$, 
\begin{eqnarray}
&&
\mbox{$\cI [f(W)|(t, Y(t))], t \geq 0$ 
is a continuous-time martingale, and}
\nonumber\\
\label{eqn:int4}
&& \cI [f(W) |(0, Y(0))]
=f(Y(0)).
\end{eqnarray}
Note that by setting $f \equiv 1$ in (\ref{eqn:int4}) we have
$\cI[1|(t, Y(t))] \equiv 1, t \geq 0$.

\begin{lem}
\label{thm:I_BM}
Let $Y$ be BM on $\R$ and $m_n, n \in \N_0$
be its associated FMPs given by Lemma \ref{thm:Hermite1}.
Then the integral transformation (\ref{eqn:int1}) with
$c=i$ and (\ref{eqn:q}) satisfies (\ref{eqn:int2}).
\end{lem}
\noindent{\it Proof.} \,
Assume $t>0$.
The Fourier transform of (\ref{eqn:G_alpha}), 
$
(1/2\pi) \int_{\R} d \alpha 
e^{-i \alpha w} G_{\alpha}(t, x), 
$
is given by (\ref{eqn:q}).
Owing to the factor $e^{-w^2/2t}$ in $q(\cdot, w|\cdot)$,
the following calculations are justified,
\begin{eqnarray}
G_{\alpha}(t,x)
&=& \sum_{n=0}^{\infty} m_{n}(t,x) \frac{\alpha^n}{n!}
\nonumber\\
&=& \int_{\R} dw \, e^{i \alpha w} q(t,w|x)
= \sum_{n=0}^{\infty} \int_{\R} dw \, 
(iw)^n \frac{\alpha^n}{n!} q(t,w|x),
\nonumber
\end{eqnarray}
which proves 
$$
m_n(t,x)=\int_{\R} dw \, (iw)^n
\frac{e^{-(ix+w)^2/2t}}{\sqrt{2 \pi t}}, \quad n \in \N_0, \quad t>0.
$$
We change the integral variable by $\widetilde{y}=ix+w$.
Since the integrand is an entire function of $\widetilde{y}$,
the obtained integral on $\widetilde{y} \in (-\infty+ix, \infty+ix)$
is equal to that on $\widetilde{y} \in \R$
by Cauchy's integral theorem,
$$
m_n(t,x)=\int_{\R} d \widetilde{y} \, (x+i \widetilde{y})^n
p_{\rm BM}(t, \widetilde{y}|0), \quad n \in \N_0, \quad t>0.
$$ 
Then the limit $t \downarrow 0$ gives
$m_n(0,x)=\int_{\R} d \widetilde{y} \, (x+i \widetilde{y})^n 
\delta_0(\{\widetilde{y}\})=x^n, n \in \N_0$
as required.
The proof is completed. \qed
\vskip 0.3cm

\begin{lem}
\label{thm:I_BESQ}
Let $Y$ be BESQ$^{(\nu)}$ on $\R_+$,
$\nu > -1$, and $m_n, n \in \N_0$
be its associated FMPs
given by Lemma \ref{thm:Laguerre1}.
If we put $c=-1$ and 
$
q(t,w|x)=q^{(\nu)}(t,w|t)
$
in (\ref{eqn:int1}), where $q^{(\nu)}$ is
given by (\ref{eqn:q_nu}),
then (\ref{eqn:int2}) is satisfied.
\end{lem}
\noindent{\it Proof.} \,
Assume $t>0$ and $x>0$.
The integral formula (\ref{eqn:int_form2})
gives a Laplace transform of (\ref{eqn:Gnu1}), 
$
G_{\alpha}^{(\nu)}(t,x)
=\int_{\R_+} dw \,
q^{(\nu)}(t,w|x) e^{-\alpha w}
$
with (\ref{eqn:q_nu}).
Owing to the factor $e^{-w/2t}$ in $q^{(\nu)}(\cdot, w|\cdot)$,
expansion with respect to $\alpha$ and integral over $w$
is exchangeable and we obtain
$$
m_n(t,x)=\int_{\R_+} dw \, (-w)^n
q^{(\nu)}(t,w|x), \quad
n \in \N_0, \quad t>0, \quad x>0.
$$
For $t>0$, (\ref{eqn:q_nu}) gives
$$
\lim_{x \downarrow 0} q^{(\nu)}(t,w|x)
=\frac{w^{\nu}e^{-w/2t}}{(2t)^{\nu+1} \Gamma(\nu+1)}
\equiv q^{(\nu)}(t,w|0),
\quad w \in \R_+.
$$
Note that $q^{(\nu)}(t,w|0)=p^{(\nu)}(t, w|0), t>0, w \in \R_+$.
Since $J_{\nu}(z)/z^{\nu}$ is an entire function of $z$,
$q^{(\nu)}(\cdot, \cdot|x)$ can be analytically continued to $x<0$;
by (\ref{eqn:I_J}) we see
$q^{(\nu)}(t, w|-|x|)=p^{(\nu)}(t, w||x|)$ with (\ref{eqn:p_nu}).
It implies $\lim_{t \downarrow 0} q^{(\nu)}(t, w|-|x|)
=p^{(\nu)}(0, w||x|)=\delta_{|x|}(\{w\})$, and we will define
$q^{(\nu)}(0, w|x)=\lim_{t \downarrow 0}q^{(\nu)}(t, w|x)
=\delta_{-x}(\{w\}), x \in \R$.
It gives $\lim_{t \downarrow 0} m_n(t,x)=x^n, n \in \N_0$
as required. Then the proof is completed.
\qed
\vskip 0.3cm

Now we can prove the following.
\begin{prop}
\label{thm:cM1}
For Processes 1 and 2, if $\xi \in \mM_0(S)$, 
$\cM_{\xi}$ given by (\ref{eqn:cMs})
satisfies the conditions (M1)-(M3).
\end{prop}
\noindent{\it Proof.} \,
For Process 1 and 2, for each $1 \leq k \leq N$,
$\Phi_{\xi}^{u_k}(z)$ is a polynomial of $z$ with degree $N-1$.
Then it can be written as
$\Phi_{\xi}^{u_k}(z)=\sum_{n=0}^{N-1} a_n z^n$,
where the coefficients $a_n \in \C$ are functions of
$\xi$ ({\it i.e.}, $u_j$'s).
By definition (\ref{eqn:cMs}) and 
Lemmas \ref{thm:I_BM} and \ref{thm:I_BESQ}, 
$\cM_{\xi}^{u_k}(t, x)=\sum_{n=0}^{N-1} a_n m_n(t, x)$,
where $m_n(t, x)$ are FMPs
given by Lemmas \ref{thm:Hermite1} and \ref{thm:Laguerre1}.
Then the condition (M1) is proved.
Since $\xi \in \mM_0(S)$ is assumed,
the set of zeros of $\Phi_{\xi}^{u_k}(z)$ 
is different from that of $\Phi_{\xi}^{u_{\ell}}(z)$,
if $k \not= \ell$, and then the condition (M2) is satisfied.
By (\ref{eqn:Phi}), $\Phi_{\xi}^{u_k}(u_j)=\delta_{jk}, 1 \leq j, k \leq N$.
Then the condition (M3) is also satisfied.
\qed
\vskip 0.3cm

The determinantal martingale (\ref{eqn:D1})
is now given by
\begin{equation}
\cD_{\xi}(t, \Y(t))
= \det_{1 \leq j, k \leq N}
\Big[ \cI[ \Phi_{\xi}^{u_k}(W)|
(t, Y_j(t))] \Big],
\quad t \in [0, \infty).
\label{eqn:D3}
\end{equation}
The integral transformation (\ref{eqn:int1}) is extended
to the linear integral transformation of functions
of $\x \in S^N$ so that, if
$F^{(k)}(\x)=\prod_{j=1}^N f_j^{(k)}(x_j)$, $k=1,2$,
then
$$
\cI \left[ F^{(k)}(\bW) \left| \{(t_{\ell}, x_{\ell})\}_{\ell=1}^N  \right. \right]
= \prod_{j=1}^N \cI \left[ \left. f^{(k)}_j(W_j) \right| (t_{j}, x_{j}) \right],
\quad k=1,2,
$$
and
\begin{eqnarray}
&& \cI \Big[c_1 F^{(1)}(\bW) +c_2 F^{(2)}(\bW) \left| \{(t_{\ell}, x_{\ell})\}_{\ell=1}^N 
\right. \Big]
\nonumber\\
&& \quad
= c_1 \cI \Big[F^{(1)}(\bW) \left| \{(t_{\ell}, 
x_{\ell})\}_{\ell=1}^N \right. \Big]
+ c_2 \cI \Big[F^{(2)}(\bW) \left| \{(t_{\ell}, 
x_{\ell})\}_{\ell=1}^N \right. \Big],
\nonumber
\end{eqnarray}
$c_1, c_2 \in \C$, 
for $0 < t_j < \infty, 1 \leq j \leq N$,
where $\bW=(W_1, \dots, W_N) \in S^N$.
In particular, if $t_{\ell}=t, 1 \leq {^{\forall}\ell} \leq N$, we write
$\cI[\cdot | \{(t_{\ell}, x_{\ell})\}_{\ell=1}^N]$ simply as
$\cI[\cdot|(t, \x)]$ with $\x=(x_1, \dots, x_N)$.
Then, by the multilinearity of determinant,
(\ref{eqn:D3}) is written as
\begin{equation}
\cD_{\xi}(t, \Y(t))
= \cI \left[ \left. \det_{1 \leq j, k \leq N}
\left[ \Phi_{\xi}^{u_k}(W_j) \right] \right|
(t, \Y(t))] \right],
\quad t \in [0, \infty).
\label{eqn:D4}
\end{equation}

\subsection{Absorbing Brownian motion in Weyl alcove
and martingale \label{sec:alcove}}

Consider the following bounded region in $\R^N$,
\begin{equation}
\cA^{A_{N-1}}_{2 \pi r}=\{ \x=(x_1, \dots, x_N) \in \R^N :
x_1< x_2 < \cdots < x_N < x_1+2 \pi r \}.
\label{eqn:WeylA}
\end{equation}
It is called a scaled {\it alcove} of the affine Weyl group of type
$A_{N-1}$ (with scale $2 \pi r$) \cite{Gra02,Kra07}.
The $N$-particle system of Brownian motions
on $\rS^1(r)$, whose process is killed when any pair of 
particles collide, is identified with the
$N$-dimensional absorbing Brownian motion in 
the Weyl alcove $\cA^{A_{N-1}}_{2 \pi r}$.
Let $q_N^r(t, \y|\x)$ be the TD of this $N$-particle process.
For $\x, \y \in \cA^{A_{N-1}}_{2 \pi r}$, it solves the $N$-dimensional
diffusion equation
\begin{equation}
\frac{\partial q_N^r(t, \y|\x)}{\partial t}
=\frac{1}{2} \sum_{j=1}^N \frac{\partial^2 q_N^r(t, \y|\x)}{\partial x_j^2},
\quad t \in [0, \infty),
\label{eqn:Ndiffeq}
\end{equation}
under the boundary condition
$$
q_N^r(t, \y|\x)=0,
\quad \x \in \partial \cA^{A_{N-1}}_{2 \pi r},
\quad \y \in \cA^{A_{N-1}}_{2 \pi r},
\quad t \in [0, \infty), 
$$
and the initial condition
\begin{equation}
q_N^r(0, \y|\x)=\prod_{j=1}^N \delta_{x_j}(\{y_j\}).
\label{eqn:initial_c1}
\end{equation}
When $N$ is odd, the Karlin-McGregor argument \cite{KM59}
immediately gives the determinantal expression to $q_N^r$ as
\begin{equation}
q_N^r(t, \y|\x)=
\det_{1 \leq j, k \leq N} 
\left[ \sum_{\ell \in \Z} p_{\rm BM}(t, y_j+2 \pi r \ell|x_k) \right],
\label{eqn:KM1}
\end{equation}
where the $(j, k)$-entry of determinant is 
obtained by wrapping the TDs of BM given by
(\ref{eqn:p}) \cite{BPY01}. 
When $N$ is even, however, (\ref{eqn:KM1}) does not hold
as pointed by Forrester \cite{For90a}.
If we moderate the initial condition (\ref{eqn:initial_c1}) as
\begin{equation}
q_N^r(0, \y|\x)=\sum_{\sigma \in \cS_N} \prod_{j=1}^N \delta_{x_j}(\{y_{\sigma(j)}\}),
\label{eqn:initial_c2}
\end{equation}
that is, we do not require that for $1 \leq j \leq N$ the particle
started at $x_j$ ends at $y_j$, but only that
it ends at $y_k$ for some $1 \leq k \leq N$, 
then we can obtain the determinantal formula,
\begin{equation}
q_N^r(t, \y|\x)=
\det_{1 \leq j, k \leq N} 
\left[ \sum_{\ell \in \Z} (-1)^{\ell}
p_{\rm BM}(t, y_j+2 \pi r \ell|x_k) \right],
\label{eqn:LW01}
\end{equation}
for even $N$, $\x, \y \in \cA^{A_{N-1}}_{2 \pi r}$,
$t \in [0, \infty)$.
A discrete analogue of this problem was first solved
by Fulmek \cite{Ful04} and as taking the continuum limit
Liechty and Wang gave
proof of (\ref{eqn:LW01}) 
 (Proposition 1.1 in \cite{LW13}).
Then, if we define $p^r(\cdot, \cdot|\cdot; N)$ as (\ref{eqn:prB1}),
the results are summarized as follows,
\begin{equation}
q_N^r(t, \y|\x)
=\det_{1 \leq j, k \leq N}
[ p^r(t, y_j|x_k; N)],
\quad \x, \y \in \cA^{A_{N-1}}_{2 \pi r},
\quad t \in [0, \infty).
\label{eqn:qNr1}
\end{equation}

We give the following remarks on $p^r$.
Let $v, \tau \in \C$ and put
$$
z=z(v)=e^{\pi i v}, \quad q=q(\tau)=e^{\pi i \tau}.
$$
The Jacobi theta functions $\vartheta_{\mu}, \mu=0,1,2,3$ are defined as
\begin{eqnarray}
\vartheta_0(v; \tau)
&=& \sum_{n \in \Z} (-1)^n q^{n^2} z^{2 n},
\nonumber\\
\vartheta_1(v; \tau)
&=& -i q^{1/4} z \vartheta_0 \left( v+\frac{\tau}{2}; \tau \right)
= i \sum_{n \in \Z} (-1)^n q^{(n-(1/2))^2} z^{2n-1},
\nonumber\\
\vartheta_2(v; \tau)
&=& q^{1/4} z \vartheta_0 \left( v+\frac{1+\tau}{2}; \tau \right)
= \sum_{n \in \Z} q^{(n-(1/2))^2} z^{2n-1},
\nonumber\\
\vartheta_3(v; \tau)
&=& \vartheta_0 \left( v+\frac{1}{2}; \tau \right)
= \sum_{n \in \Z} q^{n^2} z^{2n}.
\nonumber
\end{eqnarray}
(Note that the present functions $\vartheta_{\mu}(v; \tau)$ are represented
as $\vartheta_{\mu}(\pi v,q)$, $\mu=0,1,2,3$ in \cite{WW27}.)
For $\Im \tau >0$, $\vartheta_{\mu}(v; \tau)$
are holomorphic for $|v| < \infty$
and solve the partial differential equation
$$
\frac{\partial \vartheta_{\mu}(v; \tau)}{\partial \tau}=
\frac{1}{4 \pi i} \frac{\partial^2 \vartheta_{\mu}(v; \tau)}{\partial v^2}.
$$
By definition, we can readily confirm that
$$
p^r(t, y|x; N)= \left\{
\begin{array}{ll}
\displaystyle{
p_{\rm BM}(t, y|x) \vartheta_3 \left( \frac{i(y-x)r}{t}; \frac{2 \pi i r^2}{t} \right)
}, \quad & \mbox{if $N$ is odd}, \cr
\displaystyle{
p_{\rm BM}(t, y|x) \vartheta_0 \left( \frac{i(y-x)r}{t}; \frac{2 \pi i r^2}{t} \right)
}, \quad & \mbox{if $N$ is even}.
\end{array} \right.
$$
The equalities
\begin{eqnarray}
\vartheta_0(v; \tau)
&=& e^{\pi i/4} \tau^{-1/2} e^{-\pi i v^2/\tau}
\vartheta_2 \left( \frac{v}{\tau}; - \frac{1}{\tau} \right),
\nonumber\\
\vartheta_3(v; \tau)
&=& e^{\pi i/4} \tau^{-1/2} e^{-\pi i v^2/\tau}
\vartheta_3 \left( \frac{v}{\tau}; - \frac{1}{\tau} \right),
\nonumber
\end{eqnarray}
are established and they are called
{\it Jacobi's imaginary transformations} \cite{WW27}.
Then we obtain
\begin{equation}
p^r(t, y|x; N)= \left\{
\begin{array}{ll}
\displaystyle{
\frac{1}{2 \pi r} \vartheta_3 \left( \frac{y-x}{2 \pi r}; \frac{it}{2 \pi r^2} \right)
}, \quad & \mbox{if $N$ is odd}, \cr
\displaystyle{
\frac{1}{2 \pi r} \vartheta_2 \left( \frac{y-x}{2 \pi r}; \frac{it}{2 \pi r^2} \right)
}, \quad & \mbox{if $N$ is even}.
\end{array} \right.
\label{eqn:prB3}
\end{equation}
If we introduce the notation
\begin{equation}
\sigma_N(m)= \left\{
\begin{array}{ll}
m, \quad & \mbox{when $N$ is odd}, \cr
m-1/2, \quad & \mbox{when $N$ is even},
\end{array} \right.
\label{eqn:spin}
\end{equation}
for $m \in \Z$, then (\ref{eqn:prB3}) is written as
\begin{equation}
p^r(t, y|x; N)
= \frac{1}{2 \pi r}
\sum_{\ell \in \Z} e^{-\sigma_N(\ell)^2 t/2r^2+i \sigma_N(\ell) (y-x)/r}.
\label{eqn:pR2}
\end{equation}
It should be noted that this expression is
found in Nagao and Forrester \cite{NF02}.

Let $Y(t), t \in [0, \infty)$ be the Markov process
in $S=[0, 2 \pi r)$ such that the TD is given by (\ref{eqn:prB1}).
We have called it `BM on $[0, 2 \pi r)$'.
We prove the following.
\begin{prop}
\label{thm:cM2}
For Process 3, if $\xi \in \mM_0([0, 2 \pi r))$, 
$\cM_{\xi}$ given by (\ref{eqn:cMs})
satisfies the conditions (M1)-(M3).
\end{prop}
\noindent{\it Proof.} \,
For Process 3, (\ref{eqn:CBM}) with (\ref{eqn:Phi}) gives
\begin{equation}
\cM_{\xi}^{u_k}(t, y)
= \widetilde{\rE} \left[
\prod_{\substack{1 \leq \ell \leq N, \cr \ell \not= k}}
\frac{\sin((y+i \widetilde{Y}(t)-u_{\ell})/2r)}
{\sin((u_k-u_{\ell})/2r)} \right],
\quad 1 \leq k \leq N.
\label{eqn:Msin1}
\end{equation}
It implies that for $\ell \in \Z$,
\begin{eqnarray}
\cM_{\xi}^{u_k}(t, y+2 \pi r \ell)
&=& (-1)^{\ell (N-1)} \cM_{\xi}^{u_k}(t, y)
\nonumber\\
&=& 
\left\{
\begin{array}{ll}
\cM_{\xi}^{u_k}(t, y), \quad & \mbox{if $N$ is odd}, \cr
(-1)^{\ell} \cM_{\xi}^{u_k}(t, y), \quad & \mbox{if $N$ is even}.
\end{array} \right.
\label{eqn:Mparity}
\end{eqnarray}
Then, for $0 \leq s \leq t, 1 \leq k \leq N$,
(\ref{eqn:prB1}) gives
\begin{eqnarray}
&& \rE[\cM_{\xi}^{u_k}(t, Y(t)) | \cF(s)]
= \int_{0}^{2 \pi r} dy \,
\cM_{\xi}^{u_k}(t, y) p^r(t-s, y|Y(s); N)
\nonumber\\
&& =
\left\{
\begin{array}{ll}
\displaystyle{
\sum_{\ell \in \Z} \int_{2 \pi r \ell}^{2 \pi r (\ell+1)} dy \,
\cM_{\xi}^{u_k}(t, y-2 \pi r \ell) p_{\rm BM}(t-s, y|Y(s))
}, \quad & \mbox{if $N$ is odd}, \cr
\displaystyle{
\sum_{\ell \in \Z} \int_{2 \pi r \ell}^{2 \pi r (\ell+1)} dy \, (-1)^{\ell}
\cM_{\xi}^{u_k}(t, y-2 \pi r \ell) p_{\rm BM}(t-s, y|Y(s))
}, \quad & \mbox{if $N$ is even}.
\end{array} \right.
\nonumber
\end{eqnarray}
By (\ref{eqn:Mparity}), it is equal to
$$
\int_{\R} dy \, \cM_{\xi}^{u_k}(t, y) 
p_{\rm BM}(t-s, y|Y(s)) \quad
\mbox{a.s.}
$$
We can expand the function $\Phi_{\xi}^{u_k}(y+i \widetilde{y})$
for Process 3 as
$$
\prod_{\substack{1 \leq \ell \leq N, \cr \ell \not= k}}
\frac{\sin((y+i \widetilde{y}-u_{\ell})/2r)}{\sin((u_k-u_{\ell})/2r)}
=\sum_{n=-(N-1)}^{N-1} b_n e^{i n(y+i \widetilde{y})/2r},
$$
where the coefficients $b_n$ are functions
of $u_j$'s. Then (\ref{eqn:Msin1}) gives
\begin{eqnarray}
\cM_{\xi}^{u_k}(t, y)
&=& \sum_{n=-(N-1)}^{N-1} b_n e^{i n y/2r} \widetilde{\rE}
[e^{-n \widetilde{Y}(t)/2r} ] 
\nonumber\\
&=& \sum_{n=-(N-1)}^{N-1} b_n G_{i n / 2 r}(t, y),
\nonumber
\end{eqnarray}
where $G_{\alpha}$ is defined by (\ref{eqn:G_alpha}).
As shown by (\ref{eqn:G_martingale}) for BM on $\R$, $Y(t), t \in [0, \infty)$ , 
$$
\int_{\R} dy \, G_{\alpha}(t, y) p_{\rm BM}(t-s, y|x)
= G_{\alpha}(s, x), \quad 0 \leq s \leq t,
$$
for any $\alpha \in \C$. Then we can conclude
$$
\rE[\cM_{\xi}^{u_k}(t, Y(t)) | \cF(s)]
= \cM_{\xi}^{u_k}(s, Y(s)) \quad \mbox{a.s.}
$$
for $0 \leq s \leq t$, that is, (M1) is proved.
Since $\xi \in \mM_0([0, 2 \pi r))$ is assumed,
the set of zeros of $\Phi_{\xi}^{u_k}(z)$ 
is different from that of $\Phi_{\xi}^{u_{\ell}}(z)$,
if $k \not= \ell$, and the condition (M2) is satisfied.
By (\ref{eqn:Phi}), $\Phi_{\xi}^{u_k}(u_j)=\delta_{jk}, 1 \leq j, k \leq N$.
Then the condition (M3) is also satisfied.
\qed
\vskip 0.3cm

\SSC{Proof of Theorem \ref{thm:3_process} \label{sec:proof2}}

\subsection{Applications of Krattenthaler's determinant identity}

In the present paper
the Vandermonde determinant is denoted by
\begin{equation}
h(\x)=\det_{1 \leq j, k \leq N}
[x_j^{k-1}] = \prod_{1 \leq j< k \leq N} (x_k-x_j),
\label{eqn:Vand}
\end{equation}
for $\x=(x_1, \dots, x_N) \in \R^N$.
We define the degree of a Laurent polynomial
$\pi(x)=\sum_{j=M}^N a_j x^j$,
$M, N \in \N$, $M<N$, $a_j \in \R$, $a_N \not=0$
to be deg $\pi=N$.
The following determinant identity is proved 
in \cite{Kra99} (Lemma 5).
Here we use the convention that empty products equal 1.

\begin{lem}
\label{thm:Kra1}
Let $x_1, x_2, \dots, x_N$, $a_2, a_3, \dots, a_N$, $c$
be indeterminates.
If $\pi_0, \pi_1, \dots, \pi_{N-1}$ are Laurent polynomials
with deg $\pi_k \leq k$ and
$\pi_k(c/x)=\pi_k(x)$ for $k=0,1, \dots, N-1$, then
\begin{eqnarray}
&& \det_{1 \leq j,k \leq N}
\Big[ (x_j+a_N)(x_j+a_{N-1}) \cdots (x_j+a_{k+1})
\nonumber\\
&& \qquad \qquad \times (c/x_j+a_N)(c/x_j+a_{N-1}) 
\cdots (c/x_j+a_{k+1}) \cdot \pi_{k-1}(x_j) \Big]
\nonumber\\
&& \qquad 
= \prod_{1 \leq j < k \leq N}
(x_j-x_k)(1-c/(x_j x_k))
\prod_{j=1}^N a_j^{j-1} \prod_{j=1}^N \pi_{j-1}(-a_j).
\label{eqn:Kra1}
\end{eqnarray}
\end{lem}
\noindent

As a corollary, we obtain the following identities.

\begin{cor}
\label{thm:det_iden}
For $\x \in \R^N, \u \in \W_N(\R)$,
\begin{equation}
\det_{1 \leq j, k \leq N}
\left[ 
\prod_{1 \leq \ell \leq N, \ell \not=k}
\frac{x_j-u_{\ell}}{u_k-u_{\ell}}
 \right]
= \frac{h(\x)}{h(\u)},
\label{eqn:det_iden1}
\end{equation}
and for $\y \in [0, 2 \pi r)^N, \w \in \W_N([0, 2 \pi r))$, $r >0$, 
\begin{equation}
\det_{1 \leq j, k \leq N}
\left[ 
\prod_{1 \leq \ell \leq N, \ell \not= k}
\frac{\sin((y_j-w_{\ell})/2r)}{\sin((w_k-w_{\ell})/2r)} \right]
=\prod_{1 \leq j < k \leq N}
\frac{\sin((y_k-y_j)/2r)}{\sin((w_k-w_j)/2r)}.
\label{eqn:det_iden2}
\end{equation}
\end{cor}
\noindent{\it Proof.} \,
First we put
$a_k=-u_k, 2 \leq k \leq N$, $c=0$, and
$$
\pi_k(x)=\prod_{\ell=1}^k(x-u_{\ell}), \quad 0 \leq k \leq N-1
$$ in
(\ref{eqn:Kra1}). Then we have
$$
\det_{1 \leq j, k \leq N}
\left[ 
\prod_{1 \leq \ell \leq N, \ell \not=k}
(x_j-u_{\ell}) \right]
 =(-1)^{N(N-1)/2} h(\x) h(\u).
$$
If we divide the both sides by
$\prod_{1 \leq j, k \leq N, j \not= k}(u_k-u_j)$, 
then (\ref{eqn:det_iden1}) is obtained \cite{KT13}.

Next we put
$a_k=-u_k, 2 \leq k \leq N$, $c=1$, and
$$
\pi_k(x)=\prod_{\ell=1}^k(x-u_{\ell})(1/x-u_{\ell}), \quad 0 \leq k \leq N-1
$$
in (\ref{eqn:Kra1}).
It is easy to confirm that $\pi_k(x)=\pi_k(1/x)$ is satisfied.
Then by dividing both side by
$
\prod_{1 \leq j, k \leq N, j \not= k}
(u_j-u_k)(1/u_j-u_k),
$
we have the equality
\begin{equation}
\det_{1 \leq j, k \leq N}
\left[ 
\prod_{1 \leq \ell \leq N, \ell \not= k}
\frac{(x_j-u_{\ell})(1/x_j-u_{\ell})}{(u_k-u_{\ell})(1/u_k-u_{\ell})} \right]
=\prod_{1 \leq j < k \leq N}
\frac{(x_k-x_j)(1/(x_j x_k)-1)}{(u_k-u_j)(1/(u_j u_k)-1)}.
\label{eqn:det_iden2B}
\end{equation}
Put
$$
x_j=i \tan\left( \frac{y_j}{4 r} \right),
\quad
u_j=i \tan \left( \frac{w_j}{4 r} \right),
\quad 1 \leq j \leq N
$$
in (\ref{eqn:det_iden2B}). Then we obtain (\ref{eqn:det_iden2})
and the proof is completed. \qed
\vskip 0.3cm

As applications of the above determinant identities,
we obtain the following expressions of
determinantal martingales (\ref{eqn:D4}).

\begin{prop}
\label{thm:DM_h}
For Processes 1 and 2, when $\y \in S^N$, $\u \in \W_N(S)$, 
\begin{equation}
\cD_{\xi}(t, \y)
= \frac{h(\y)}{h(\u)},
\quad 0 \leq t < \infty.
\label{eqn:D6}
\end{equation}
\end{prop}
\noindent{\it Proof.} \,
Using the determinant identity (\ref{eqn:det_iden1})
with (\ref{eqn:Vand}), (\ref{eqn:D4}) gives
\begin{eqnarray}
\cD_{\xi}(t, \y) 
&=& \cI \left[ \left. 
\frac{h(\bW)}{h(\u)} \right| (t, \y) \right]
\nonumber\\
&=& \cI \left[ \left.
\frac{1}{h(\u)} \det_{1 \leq j, k \leq N}
[(W_j)^{k-1}] \right| (t, \y) \right]
\nonumber\\
&=& \frac{1}{h(\u)}
\det_{1 \leq j, k \leq N}
\Big[ \cI [(W_j)^{k-1} | (t, y_j)] \Big]
\nonumber\\
&=& 
\frac{1}{h(\u)} \det_{1 \leq j, k \leq N}
[m_{k-1}(t, y_j)],
\nonumber
\end{eqnarray}
where Lemmas \ref{thm:I_BM} and \ref{thm:I_BESQ} were used.
By multilinearity of determinant, the Vandermonde determinant 
$\det_{1 \leq j, k \leq N} [y_j^{k-1} ] $ does not change
by replacing $y_j^{k-1}$ by any monic polynomial of $y_j$ of 
degree $k-1$, $1 \leq j, k \leq N$.
Since $m_{k-1}(t,y_j)$ is a monic polynomial of $y_j$
of degree $k-1$, the above is equal to
$\det_{1 \leq j, k \leq N}[y_j^{k-1}]/h(\u)$.
Then the proof is completed. \qed
\vskip 0.3cm

\begin{prop}
\label{thm:DM_ht}
For $\x \in [0, 2 \pi r)^N, r>0, 0 \leq t < \infty$, 
let
\begin{equation}
h^r(t, \x)=e^{t N(N^2-1)/24 r^2}
\prod_{1 \leq j < k \leq N}
\sin \left( \frac{x_k-x_j}{2r} \right).
\label{eqn:ht}
\end{equation}
Then for Process 3, when $\y \in [0, 2 \pi r)^N$, $\u \in \W_N([0, 2 \pi r))$,
\begin{equation}
\cD_{\xi}(t, \y)
= \frac{h^r(t,\y)}{h^r(0,\u)},
\quad 0 \leq t < \infty.
\label{eqn:Dt}
\end{equation}
\end{prop}
\noindent{\it Proof.} \,
By definition (\ref{eqn:Phi}) of $\Phi_{\xi}^v$ for Process 3
and (\ref{eqn:det_iden2}) of Corollary \ref{thm:det_iden},
$$
\cD_{\xi}(t, \y)
=\cI \left[ \left. 
\det_{1 \leq j, k \leq N}
[\Phi_{\xi}^{u_k}(W_j)] \right| (t, \y) \right]
=\frac{\cI[h^r(0, \bW) | (t, \y)]}{h^r(0, \u)}. 
$$
Since
\begin{equation}
h^r(0, \x)=(2i)^{-N(N-1)/2}
\det_{1 \leq j, k \leq N} [e^{i(2k-N-1)x_j/2r}],
\label{eqn:Vand0}
\end{equation}
multilinearity of determinant and the integral transformation
$\cI[\cdot|(t, \y)]$ gives
\begin{equation}
\cI[h^r(0, \bW) | (t, \y)]
= (2i)^{-N(N-1)/2}
\det_{1 \leq j, k \leq N} \Bigg[
\cI[e^{i(2k-N-1) W_j/2r} | (t, y_j)]
\Bigg].
\label{eqn:Dt3}
\end{equation}
By (\ref{eqn:I}) with (\ref{eqn:q}) for Process 3, we see
$$
\cI[e^{i(2k-N-1) W_j/2r} | (t, y_j)]
=G_{i(2k-N-1)/2r}(t, y_j)
$$ 
with (\ref{eqn:G_alpha}),
and we find that
(\ref{eqn:Dt3}) is equal to
\begin{eqnarray}
&& (2 i)^{-N(N-1)/2}
\det_{1 \leq j, k \leq N}
\Big[ e^{i(2k-N-1) y_j/2r}
e^{t(2k-N-1)^2/8r^2} \Big]
\nonumber\\
&& \qquad 
= \prod_{k=1}^N e^{t(2k-N-1)^2/8r^2}
(2i)^{-N(N-1)/2}
\det_{1 \leq j, k \leq N}
\Big[ e^{i(2k-N-1) y_j/2r} \Big]
\nonumber\\
&& \qquad
= e^{ t \sum_{k=1}^N (2k-N-1)^2/8r^2}
h^r(0, \y),
\nonumber
\end{eqnarray}
where (\ref{eqn:Vand0}) was again used.
For 
$
\sum_{k=1}^N (2k-N-1)^2
=N(N^2-1)/3,
$
the proof is completed. \qed
\vskip 0.3cm
\noindent{\bf Remark 1} \,
The function $h^r(\cdot, \x)$ is not harmonic,
but $h^r(t, \Y(t))=h^r(0, \u) \cD_{\xi}(t, \Y(t))$ was proved to be
a martingale by Proposition \ref{thm:cM2},
when $Y_j(t), 1 \leq j \leq N$ are BM on $[0, 2 \pi r)$.
Then It\^o's formula implies
\begin{equation}
\left( \frac{\partial}{\partial t} + \frac{1}{2} \Delta \right) h^r(t, \x)=0,
\label{eqn:heq}
\end{equation}
where $\Delta=\sum_{j=1}^N \partial^2/\partial x_j^2$.
In this sense DMR for Process 3 
is a time-dependent extension of $h$-transform.

\subsection{Proof of Theorem \ref{thm:3_process}
for Processes 1 and 2}

We can prove that Process 1 (noncolliding Brownian motion)
is obtained as 
an $h$-transform of the absorbing Brownian motion in the
Weyl chamber $\W_N(\R)$ \cite{Gra99}.
Similarly, Process 2 (noncolliding BESQ$^{(\nu)}$)
is realized as an $h$-transform of the absorbing
BESQ$^{(\nu)}(t)$
in $\W_N(\R_+)$ \cite{KO01}.
In both cases, the harmonic function is
given by the Vandermonde determinant (\ref{eqn:Vand}).
For these two processes,
put
$$
\tau = \inf \{ t > 0 :
\Y(t) \notin \W_N(S) \}.
$$
Remember the note on measurable functions 
given above Definition \ref{thm:DM_rep} 
in Section \ref{sec:Introduction}.
Then for
$\xi=\sum_{j=1}^N \delta_{u_j} \in \mM_0(S)$,
$\u=(u_1, \dots, u_N) \in \W_N(S)$, 
the above statement is expressed by the equality
\begin{equation}
\E_{\xi} \left[
\prod_{m=1}^M g_{m}(\X(t_m)) \right]
=\rE_{\u} \left[
\1(\tau > t_{M})
\prod_{m=1}^{M} g_m(\Y(t_m)) 
\frac{h(\Y(t_M))}{h(\u)} \right].
\label{eqn:harmonic1}
\end{equation}

We introduce the stopping times
$$
\tau_{jk}=\inf\{t> 0 : Y_j(t)= Y_k(t)\}, \quad
1\leq j < k \leq N.
$$
Let $\pi_{jk} \in \cS_{N}$ be the permutation of $(j,k), 1 \leq j,k \leq N$.
Note that in a configuration $\u'$, if $u'_j=u'_k, j \not=k$, then
$\pi_{jk}(\u')= \u'$, 
and the processes $\Y(\cdot)$ and $\pi_{jk}(\Y(\cdot))$ are equivalent 
in probability law under the probability measure $\rP_{\u'}$.
By the strong Markov property of the process $\Y(\cdot)$
and by the fact that $h$ is anti-symmetric and $g_m, 1 \leq m \leq M$ are symmetric, 
$$
\rE_{\u} \left[
\1(\tau =\tau_{jk} < t_{M})
\prod_{m=1}^{M} g_m(\Y(t_m)) 
\frac{h(\Y(t_M))}{h(\u)} \right]=0,
\quad 1 \leq j < k \leq N.
$$
Since $\rP_{\u}(\tau_{jk}= \tau_{j' k'})=0$ if $(j,k)\not=(j',k')$,
and 
$
\tau= \min_{1\leq j<k \leq N} \tau_{jk},
$
$$
\rE_{\u} \left[
\1(\tau < t_{M})
\prod_{m=1}^{M} g_m(\Y(t_m)) 
\frac{h(\Y(t_M))}{h(\u)} \right]=0.
$$
Hence, (\ref{eqn:harmonic1}) equals
$$
\rE_{\u} \left[
\prod_{m=1}^{M} g_m(\Y(t_m)) 
\frac{h(\Y(t_M))}{h(\u)} \right].
$$
By the equality (\ref{eqn:D6})
of Proposition \ref{thm:DM_h}, 
we obtain a required DMR.
The proof is then completed. \qed
\vskip 0.3cm

\subsection{Proof of Theorem \ref{thm:3_process}
for Process 3}

The backward Kolmogorov equation for 
(\ref{eqn:noncollBM_S1}) is given as
\begin{equation}
-\frac{\partial u(s, \x)}{\partial s}= 
\frac{1}{2} \sum_{j=1}^N \frac{\partial^2 u(s,\x)}{\partial x_j^2} 
+ \frac{1}{2 r} \sum_{\substack{1 \leq j, k \leq N, \cr j \not= k}}
\cot \left( \frac{x_j-x_k}{2r} \right)
\frac{\partial u(s,\x)}{\partial x_j}.
\label{eqn:Keq1}
\end{equation}
By definition (\ref{eqn:noncollBM_S1b}),
the TD of Process 3, denoted by
$p_N^r(t-s, \y|\x), t \geq s \geq 0$,
is obtained as the solution of (\ref{eqn:Keq1})
for $\x, \y \in \W_N([0, 2 \pi r))$, which
satisfies the condition
\begin{equation}
\lim_{s \uparrow t} u(s, \x)=
\sum_{\sigma \in \cS_N} \prod_{j=1}^N \delta_{x_j}(\{y_{\sigma(j)}\}).
\label{eqn:Keq2}
\end{equation}
Since we study the $\mM([0, 2 \pi r))$-valued process (\ref{eqn:Xi1}),
the usual condition 
$\lim_{s \uparrow t} u(s, \x)=\prod_{j=1}^N \delta_{x_j}(\{y_j \})$
for labeled particle system can be moderated as (\ref{eqn:Keq2}).
Let $q_N^r(t, \y|\x)$ be the TD of the absorbing Brownian motion
in the alcove $\cA^{A_{N-1}}_{2 \pi r}$ 
given by (\ref{eqn:qNr1}) with (\ref{eqn:prB1}). 

\begin{lem}
\label{thm:K_sol1}
The TD of Process 3
is given by
$$
p_N^{r}(t-s, \y|\x)
= \frac{h^r(t, \y)}{h^r(s, \x)}
q_N^{r}(t-s, \y|\x),
\quad 0 \leq s \leq t < \infty, \quad \x, \y \in \W_N([0, 2 \pi r)),
$$
where $h^r$ is given by (\ref{eqn:ht}).
\end{lem}
\noindent{\it Proof.} \,
Let
$u(s,\x)=f(s, \x) q_N^r(t-s,\y|\x)$
and put it into (\ref{eqn:Keq1})
assuming that $f$ is $\rC^1$ in $t$ and $\rC^2$ in $\x$.
Then we have
\begin{eqnarray} 
-\frac{\partial f(s,\x)}{\partial s} q_N^r(t-s, \y|\x)
&=& \frac{1}{2} q_N^r(t-s,\y|\x) \sum_{j=1}^N \frac{\partial^2 f(s,\x)}{\partial x_j^2}
+\sum_{j=1}^N \frac{\partial f(s,\x)}{\partial x_j}
\frac{\partial q_N^r(t-s,\y|\x)}{\partial x_j}
\nonumber\\
&+& \frac{1}{2r} q_N^r(t-s,\y|\x) 
\sum_{\substack{1 \leq j, k \leq N, \cr j \not= k}}
\cot \left( \frac{x_j-x_k}{2r} \right) \frac{\partial f(s,\x)}{\partial x_j}
\nonumber\\
&+& \frac{1}{2r} f(s,\x) 
\sum_{\substack{1 \leq j, k \leq N, \cr j \not= k}}
\cot \left( \frac{x_j-x_k}{2r} \right) \frac{\partial q_N^r(t-s,\y|\x)}{\partial x_j},
\label{eqn:eqB2}
\end{eqnarray}
since $q_N^r(t-s,\y|\x)$ satisfies the diffusion equation (\ref{eqn:Ndiffeq}).
Now we put
$f(s,\x)=g(s)/h^r(0,\x)$,
where $g$ is a $\rC^1$ function to be determined.
Since
$$
\frac{\partial f(s,\x)}{\partial x_j}
=-\frac{1}{2r} \sum_{\substack{1 \leq k \leq N, \cr j \not= k}}
\cot \left( \frac{x_j-x_k}{2r} \right) f(s, \x),
$$ 
(\ref{eqn:eqB2}) turns to be
\begin{eqnarray}
&& -8 r^2 \frac{d \log g(s)}{d s}
=
\sum_{j=1}^N \left[
\sum_{1 \leq \ell \leq N, \ell \not=j}
\frac{1}{\sin^2((x_j-x_{\ell})/2r)}
- \left( \sum_{1 \leq \ell \leq N, \ell \not=j}
\cot \left( \frac{x_j-x_{\ell}}{2r} \right) \right)^2 \right].
\nonumber\\
\label{eqn:eqB5}
\end{eqnarray}
We see that RHS of (\ref{eqn:eqB5}) is equal to 
\begin{eqnarray}
&& 
\sum_{\substack{1 \leq j, \ell \leq N \cr j \not= \ell}}
\left[ \frac{1}{\sin^2((x_j-x_{\ell})/2r)}
- \cot^2 \left( \frac{x_j-x_{\ell}}{2r} \right) \right]
\nonumber\\
&& \qquad \qquad 
- \sum_{\substack{1 \leq j, \ell, k \leq N \cr j \not= \ell \not= k}}
\cot \left( \frac{x_j-x_{\ell}}{2r} \right) 
\cot \left( \frac{x_j-x_k}{2r} \right)
\nonumber\\
&& = 
\sum_{\substack{1 \leq j, \ell \leq N \cr j \not= \ell}} 1
- \sum_{\substack{1 \leq j, \ell, k \cr j \not= \ell \not= k}}
\cot \left( \frac{x_j-x_{\ell}}{2r} \right) 
\cot \left(\frac{x_j-x_k}{2r} \right).
\nonumber
\end{eqnarray}
Here
$
\sum_{1 \leq j, \ell \leq N, j \not= \ell} 1
=N(N-1),
$
and the identity
$$
\cot(a-b)\cot(a-c)+\cot(b-a)\cot(b-c)+\cot(c-a)\cot(c-b)=-1
$$
gives
$$
\sum_{\substack{1 \leq j, \ell, k \leq N, \cr j \not= \ell \not= k}}
\cot \left( \frac{x_j-x_{\ell}}{2r} \right) 
\cot \left( \frac{x_j-x_k}{2r} \right)
=-\frac{1}{3}N(N-1)(N-2).
$$
Then (\ref{eqn:eqB5}) becomes
$
d \log g(s)/ds =-N(N^2-1)/24 r^2, 
$
which is solved as
$
g(s)=C e^{-s N(N^2-1)/24 r^2}
$
with a constant $C$.
Therefore, we have obtained the solution of (\ref{eqn:Keq1})
of the form
$u(s,\x)= C q_N^r(t-s,\y|\x)/h^r(s,\x)$.
For (\ref{eqn:initial_c2}), 
if and only if we set
$C =h^r(t,\y)$,
the condition (\ref{eqn:Keq2}) is satisfied.
Then the proof is completed. \qed
\vskip 0.3cm
\noindent{\it Proof of Theorem \ref{thm:3_process}
for Process 3}. \,
Let $\check{\Y}(t), t \in [0, \infty)$ be $N$-dimensional 
Brownian motion on $(\rS^1(r))^N$
started at $\u \in \W_N([0, 2 \pi r))$.
The expectation with respect to this process
is denoted by $\check{\rE}_{\u}$.
Put
$$
\tau=\inf \{ t > 0: 
\check{\Y}(t) \notin \cA^{A_{N-1}}_{2 \pi r} \}.
$$
Remember the notes on measurable functions
given above Definition \ref{thm:DM_rep}
and below (\ref{eqn:noncollBM_S1b})
in Section \ref{sec:Introduction}.
Then Lemma \ref{thm:K_sol1} implies 
$$
\E_{\xi} \left[ \prod_{m=1}^M g_m(\X(t_m)) \right]
=\check{\rE}_{\u} \left[
\1(\tau > t_{M}) \prod_{m=1}^M g_m(\check{\Y}(t_m))
\frac{h^r(t_M, \check{\Y}(t_M))}{h^r(0, \u)} \right],
$$
where $g_m, 1 \leq m \leq M$
are bounded symmetric and periodic measurable functions
on $[0, 2 \pi r)^N$.
By Markov property, it  is enough to consider the case 
$M=1$; $0 \leq t_1 \leq T < \infty$
$$
\E_{\xi} \left[ g_1(\X(t_1)) \right]
=\check{\rE}_{\u} \left[
\1(\tau > t_1) g_1(\check{\Y}(t_1))
\frac{h^r(t_1, \check{\Y}(t_1))}{h^r(0, \u)} \right].
$$
The definition of $h^r$ given by (\ref{eqn:ht})
and the assumption $\u \in \W_N([0, 2 \pi r))$ gives, 
$$
\mbox{(RHS)}
= \check{\rE}_{\u} \left[
\1(\tau > t_1) g_1(\check{\Y}(t_1))
\frac{|h^r(t_1, \check{\Y}(t_1))|}{h^r(0, \u)} \right].
$$
By the determinantal formula (\ref{eqn:qNr1})
of $q_N^r$ in Section \ref{sec:alcove},
the above is written as
\begin{equation}
\rE_{\u} \left[
\sum_{\sigma \in \cS_N} {\rm sgn}(\sigma)
g_1(\Y(t_1))
\1(\sigma(\Y(t_1)) \in \W_N([0, 2 \pi r)))
\frac{|h^r(t_1, \Y(t_1))|}{h^r(0, \u)} \right],
\label{eqn:Esgn}
\end{equation}
where $\Y(t_1)=(Y_1(t_1), \dots, Y_N(t_1))$ and 
each $Y_j$ is BM on $[0, 2 \pi r)$,
whose TD is given by (\ref{eqn:prB1}), $1 \leq j \leq N$.
Since
\begin{eqnarray}
&&
{\rm sgn}(\sigma) \1(\sigma(\Y(t_1)) \in \W_N([0, 2 \pi r)) )
|h^r(t_1, \Y(t_1))|
\nonumber\\
&& \quad 
= \1(\sigma(\Y(t_1)) \in \W_N([0, 2 \pi r)) )
h^r(t_1, \Y(t_1)), \quad \sigma \in \cS_N, 
\nonumber
\end{eqnarray}
(\ref{eqn:Esgn}) is equal to
\begin{eqnarray}
&&
\rE_{\u} \left[
\sum_{\sigma \in \cS_N} \1(\sigma(\Y(t_1)) \in \W_N([0, 2 \pi r)))
g_1(\Y(t_1)) \frac{h^r(t_1, \Y(t_1))}{h^r(0, \u)} \right]
\nonumber\\
&& \quad 
= \rE_{\u} \left[ g_1(\Y(t_1))
\frac{h^r(t_1, \Y(t_1))}{h^r(0, \u)} \right].
\nonumber
\end{eqnarray}
By Markov property, this result is readily generalized as
$$
\E_{\xi} \left[
\prod_{m=1}^M g_{m}(\X(t_m)) \right]
=\rE_{\u} \left[
\prod_{m=1}^{M} g_m(\Y(t_m)) 
\frac{h^r(t_M,\Y(t_M))}{h^r(0,\u)} \right],
$$
for an arbitrary $M \in \N,
0 \leq t_1 < \cdots < t_M \leq T < \infty$. 
By the equality (\ref{eqn:Dt})
of Proposition \ref{thm:DM_ht}, 
the theorem is concluded. \qed

\SSC{Martingales for Configurations with Multiple Points \label{sec:multi}}
For general $\xi \in \mM(S)$
with $\xi(S)=N < \infty$, define
$\supp \xi=\{x \in S : \xi(x) > 0\}$ and let
$
\xi_{*}(\cdot)=\sum_{v \in \supp \xi} \delta_{v}(\cdot).
$
For $s \in [0, \infty)$, $v, x \in S$, $z, \zeta \in \C$, let
\begin{equation}
\phi_{\xi}^{v}((s,x);z,\zeta)
=\left\{ \begin{array}{ll}
\displaystyle{
\frac{p(s,x|\zeta)}{p(s,x|v)} \frac{1}{z-\zeta}
\prod_{\ell=1}^N
\frac{z-u_{\ell}}{\zeta-u_{\ell}},
}
& \, \mbox{for Processes 1 and 2}, \cr
\displaystyle{
\frac{p(s,x|\zeta)}{p(s,x|v)} \frac{1}{z-\zeta}
\prod_{\ell=1}^N
\frac{\sin((z-u_{\ell})/2r)}{\sin((\zeta-u_{\ell})/2r),}
}
& \, \mbox{\rm for Process 3}, \cr
\end{array} \right.
\label{eqn:phi}
\end{equation}
and
\begin{eqnarray}
\Phi_{\xi}^{v}((s,x); z) &=& \frac{1}{2 \pi i}
\oint_{C(\delta_{v})} d \zeta \, 
\phi^v_{\xi}((s,x); z, \zeta)
\nonumber\\
&=& {\rm Res} \,\Big[\phi^v_{\xi}((s,x); z, \zeta); \zeta=v \Big],
\label{eqn:PhiB}
\end{eqnarray}
where
$C(\delta_{v})$ is a closed contour on the complex plane $\C$
encircling a point $v$ on $S$
once in the positive direction. Define
\begin{equation}
\cM_{\xi}^{u}((s, x)|(t,y))
=\cI \left[\Phi_{\xi}^{u}((s,x); W) \Big| (t,y) \right],
\quad (s,x), (t,y) \in [0, \infty) \times S.
\label{eqn:cMB}
\end{equation}
Then it is easy to see that (\ref{eqn:K1}) is rewritten as
\begin{equation}
\mbK_{\xi}(s,x;t,y)
=\int_{S} \xi_{*}(dv) p(s, x|v) \cM_{\xi}^{v}((s,x)|(t,y))
- \1(s>t) p(s-t,x|y),
\label{eqn:K2}
\end{equation}
$(s,x), (t,y) \in [0, \infty) \times S$.

We note that, even though the systems of SDEs (\ref{eqn:noncollBM})-(\ref{eqn:noncollBM_S1})
cannot be solved for any initial configuration with
multiple points, $\xi \in \mM(S) \setminus \mM_0(S)$,
the kernel (\ref{eqn:K2}) with (\ref{eqn:cMB})
is bounded and integrable also for $\xi \in \mM(S) \setminus \mM_0(S)$.
Therefore, spatio-temporal correlations are given
by (\ref{eqn:rho1}) for any $0 \leq t_1 < \dots < t_M < \infty, M \in \N$
and finite-dimensional distributions 
are determined.
\begin{prop}
\label{thm:entrance_law}
Also for $\xi \in \mM \setminus \mM_0$,
the determinantal processes 
with the correlation kernels (\ref{eqn:K2})
are well-defined.
They provide the entrance laws for
the processes $(\Xi(t),  t>0, \P_{\xi})$.
\end{prop}

In order to give examples
of Proposition \ref{thm:entrance_law}, 
here we study the extreme case such that
all $N$ points are concentrated on an origin,
\begin{equation}
\xi=N \delta_0 \quad
\Longleftrightarrow \quad
\xi_{*}=\delta_0 \quad
\mbox{with} \quad
\xi(\{0\})=N.
\label{eqn:xi01}
\end{equation}
We consider Processes 1 and 2.
For (\ref{eqn:xi01}), (\ref{eqn:phi}) and (\ref{eqn:PhiB})
become
\begin{eqnarray}
\phi_{N \delta_0}^0((s,x); z, \zeta)
&=& \frac{p(s,x|\zeta)}{p(s,x|0)}
\frac{1}{z-\zeta} \left( \frac{z}{\zeta} \right)^N
\nonumber\\
&=& \frac{p(s,x|\zeta)}{p(s,x|0)}
\sum_{\ell=0}^{\infty} \frac{z^{N-\ell-1}}{\zeta^{N-\ell}},
\nonumber
\end{eqnarray}
and
\begin{eqnarray}
\Phi_{N \delta_0}^0((s,x);z) &=& 
\frac{1}{p(s,x|0)} \sum_{\ell=0}^{\infty} z^{N-\ell-1}
\frac{1}{2 \pi i} \oint_{C(\delta_0)} d \zeta \,
\frac{p(s,x|\zeta)}{\zeta^{N-\ell}}
\nonumber\\
&=& 
\frac{1}{p(s,x|0)} \sum_{\ell=0}^{N-1} z^{N-\ell-1}
\frac{1}{2 \pi i} \oint_{C(\delta_0)} d \zeta \,
\frac{p(s,x|\zeta)}{\zeta^{N-\ell}},
\label{eqn:xi03}
\end{eqnarray}
since the integrands are holomorphic when $\ell \geq N$,
where we assume $\nu > -1$
and $C(\delta_0)$ is interpreted as 
$\lim_{\varepsilon \downarrow 0} C(\delta_{\varepsilon})$
for BESQ$^{(\nu)}$. 

For BM with TD (\ref{eqn:p}),
(\ref{eqn:xi03}) gives
\begin{eqnarray}
\Phi_{N \delta_0}^0((s,x);z)
&=& \sum_{\ell=0}^{N-1} z^{N-\ell-1}
\frac{1}{2 \pi i}
\oint_{C(\delta_0)} d \zeta \,
\frac{e^{x\zeta/s-\zeta^2/2s}}{\zeta^{N-\ell}}
\nonumber\\
&=& \sum_{\ell=0}^{N-1} 
\left( \frac{z}{\sqrt{2s}} \right)^{N-\ell-1}
\frac{1}{2 \pi i}
\oint_{C(\delta_0)} d \eta \,
\frac{e^{2(x/\sqrt{2s}) \eta-\eta^2}}{\eta^{N-\ell}}
\nonumber\\
&=& \sum_{\ell=0}^{N-1} 
\left( \frac{z}{\sqrt{2s}} \right)^{N-\ell-1}
\frac{1}{(N-\ell-1)!} H_{N-\ell-1} \left(\frac{x}{\sqrt{2s}} \right),
\nonumber
\end{eqnarray}
where we have used the contour integral representation
of the Hermite polynomials \cite{Sze91}
$$
H_n(x)=\frac{n!}{2 \pi i} \oint_{C(\delta_0)} d \eta
\frac{e^{2 x \eta-\eta^2}}{\eta^{n+1}},
\quad n \in \N_0, \quad x \in \R.
$$
Thus its integral transformation is calculated as
\begin{eqnarray}
&& \cI \left[ \left.
\Phi_{N \delta_0}^0((s,x);W) \right| (t, y) \right]
\nonumber\\
&& \quad 
= \sum_{\ell=0}^{N-1} \frac{1}{(N-\ell-1)!} H_{N-\ell-1}
\left( \frac{x}{\sqrt{2s}} \right)
\frac{1}{(2s)^{(N-\ell-1)/2}}
\cI[W^{N-\ell-1}|(t,y)]
\nonumber\\
&& \quad 
= \sum_{\ell=0}^{N-1} \frac{1}{(N-\ell-1)!} H_{N-\ell-1}
\left( \frac{x}{\sqrt{2s}} \right)
\frac{1}{(2s)^{(N-\ell-1)/2}} 
m_{N-\ell-1}(t,y)
\nonumber\\
&& \quad 
= \sum_{\ell=0}^{N-1} \frac{1}{(N-\ell-1)! 2^{N-\ell-1}} 
\left( \frac{t}{s} \right)^{(N-\ell-1)/2}
H_{N-\ell-1}
\left( \frac{x}{\sqrt{2s}} \right)
H_{N-\ell-1}
\left( \frac{y}{\sqrt{2t}} \right),
\nonumber
\end{eqnarray}
where we have used Lemmas \ref{thm:Hermite1} and \ref{thm:I_BM}.
Then we obtain the following,
\begin{eqnarray}
&& \cM_{N \delta_0}^0((s,x)|(t,Y(t)))
=\sum_{n=0}^{N-1} \frac{1}{n! 2^n} m_n(s,x) m_n(t,Y(t))
\nonumber\\ 
&& \qquad \qquad =
\sqrt{\pi} e^{x^2/4s+Y(t)^2/4t}
\sum_{n=0}^{N-1} 
 \left(\frac{t}{s}\right)^{n/2}
\varphi_n \left( \frac{x}{\sqrt{2s}} \right)
\varphi_n \left( \frac{Y(t)}{\sqrt{2t}} \right),
\label{eqn:Hermite3B}
\end{eqnarray}
where
$$
\varphi_n(x)=\frac{1}{\sqrt{ \sqrt{\pi} 2^n n!}}
H_n(x) e^{-x^2/2}, \quad
n \in \N, \quad x \in \R.
$$
Similarly, for BESQ$^{(\nu)}, \nu > -1$
with TD (\ref{eqn:p_nu}),
we obtain
\begin{eqnarray}
\Phi_{N \delta_0}^{0}((s,x); z)
&=& \frac{(2s)^{\nu} \Gamma(\nu+1)}{x^{\nu/2}}
\sum_{\ell=0}^{N-1} z^{N-\ell-1}
\frac{1}{2 \pi i} \oint_{C(\delta_0)} d \zeta \,
\frac{e^{-\zeta/2s}}{\zeta^{N-\ell+\nu/2}}
I_{\nu} \left( \frac{\sqrt{x \zeta}}{s}\right)
\nonumber\\
&=& \Gamma(\nu+1) \sum_{\ell=0}^{N-1}
\left( -\frac{z}{2s} \right)^{N-\ell-1}
\frac{1}{\Gamma(N-\ell+\nu)}
L^{(\nu)}_{N-\ell-1} \left( \frac{x}{2s} \right),
\nonumber
\end{eqnarray}
where we have used the contour integral representation
of the Laguerre polynomials
$$
L_n^{(\nu)}(x)
= \frac{\Gamma(n+\nu+1)}{x^{\nu/2}}
\frac{1}{2 \pi i} \oint_{C(\delta_0)}
d \eta \,
\frac{e^{\eta}}{\eta^{n+1+\nu/2}} J_{\nu}(2 \sqrt{\eta x})
$$
with the relation (\ref{eqn:I_J}) \cite{Wat44}.
By using Lemmas \ref{thm:Laguerre1} and \ref{thm:I_BESQ}, we have
\begin{eqnarray}
&& \cM_{N \delta_0}^{0}((s,x)|(t, Y(t))
= \cI^{(\nu)} \left[ \left.
\Phi_{N \delta_0}^{0}((s,x); W) \right| (t, Y(t)) \right]
\nonumber\\
&& \quad = \Gamma(\nu+1)
\sum_{n=0}^{N-1} \frac{1}{\Gamma(n+1) \Gamma(n+\nu+1) (2s)^{2n}}
m_n(s,x) m_n \left(t, Y(t) \right)
\nonumber\\
&& \quad =
\Gamma(\nu+1) 
\left(\frac{x}{2s} \right)^{-\nu/2}
\left(\frac{Y(s)}{2s} \right)^{-\nu/2}
e^{x/4s+Y(t)/4t}
\nonumber\\
&& \qquad \qquad \times
\sum_{n=0}^{N-1} \frac{\Gamma(n+1)}{\Gamma(n+\nu+1)}
\left( \frac{t}{s} \right)^n
\varphi_n^{(\nu)} \left( \frac{x}{2s} \right)
\varphi_n^{(\nu)} \left( \frac{Y(t)}{2t} \right),
\label{eqn:Laguerre3B}
\end{eqnarray}
where
$$
\varphi_n^{(\nu)}(x)=\sqrt{\frac{\Gamma(n+1)}{\Gamma(n+\nu+1)}}
x^{\nu/2} L_n^{(\nu)}(x) e^{-x/2},
\quad n \in \N_0, \quad x \in \R_+.
$$
Since $m_n$, $n \in \N_0$ are FMPs
associated with BM on $S=\R$ or BESQ$^{(\nu)}$ on $S=\R_+$, 
the processes (\ref{eqn:Hermite3B}) and (\ref{eqn:Laguerre3B})
are continuous-time martingales.
Then we see
$$
\rE \left[ \cM_{N \delta_0}^0((s,x)|(t, \Y(t))) \right]
= \rE \left[ \cM_{N \delta_0}^0((s,x)|(0, \Y(0))) \right]=1
$$
for $ (s,x) \in [0, T] \times S$, $0 \leq t < \infty$.

By the formula (\ref{eqn:K2}), we obtain the
correlation kernels for Process 1 as
\begin{eqnarray}
\mbK_{N \delta_0}^{1}(s,x;t,y) &=& p_{\rm BM}(s,x|0) \cM_{N \delta_0}^0((s,x)| (t,y))
-\1(s>t) p_{\rm BM}(s-t, x|y)
\nonumber\\
&=& \frac{e^{-x^2/4s}}{e^{-y^2/4t}}
\bK_{H}(s,x;t,y)
\nonumber
\end{eqnarray}
with
\begin{eqnarray}
\bK_{H}(s,x;t,y)
&=& \frac{1}{\sqrt{2s}}
\sum_{n=0}^{N-1} \left( \frac{t}{s} \right)^{n/2}
\varphi_n \left( \frac{x}{\sqrt{2s}} \right)
\varphi_n \left( \frac{y}{\sqrt{2t}} \right)
\nonumber\\
&-& \1(s>t) 
\frac{1}{\sqrt{2s}}
\sum_{n=0}^{\infty} \left( \frac{t}{s} \right)^{n/2}
\varphi_n \left( \frac{x}{\sqrt{2s}} \right)
\varphi_n \left( \frac{y}{\sqrt{2t}} \right),
\nonumber
\end{eqnarray}
and for Process 2
\begin{eqnarray}
\mbK^{2}_{N \delta_0}(s,x;t,y) 
&=& p^{(\nu)}(s,x|0) \cM_{N \delta_0}^{0} ((s,x)| (t,y)) 
-\1(s>t) p^{(\nu)}(s-t, x|y)
\nonumber\\
&=&
\frac{(x/2s)^{\nu/2} e^{-x/4s}}{(y/2t)^{\nu/2} e^{-y/4t}}
\bK_{L^{(\nu)}}(s,x;t,y)
\nonumber
\end{eqnarray}
with
\begin{eqnarray}
\bK_{L^{(\nu)}}(s,x;t,y) &=&
\frac{1}{2s}
\sum_{n=0}^{N-1} \left( \frac{t}{s} \right)^{n}
\varphi^{(\nu)}_n \left( \frac{x}{\sqrt{2s}} \right)
\varphi^{(\nu)}_n \left( \frac{y}{\sqrt{2t}} \right)
\nonumber\\
&-& \1(s>t) 
\frac{1}{2s}
\sum_{n=0}^{\infty} \left( \frac{t}{s} \right)^{n}
\varphi^{(\nu)}_n \left( \frac{x}{\sqrt{2s}} \right)
\varphi^{(\nu)}_n \left( \frac{y}{\sqrt{2t}} \right),
\nonumber
\end{eqnarray}
where $\bK_{H}$ and $\bK_{L^{(\nu)}}$
are known as the {\it extended Hermite and Laguerre kernels},
respectively (see, for instance, \cite{For10}).
Here we would like to emphasize the fact that
these kernels have been derived here by not following
any `orthogonal-polynomial arguments' \cite{Meh04,For10}
but by only using proper martingales
determined by the chosen initial configuration (\ref{eqn:xi01}).
The above kernels determine the finite-dimensional distributions
and specify the entrance laws
for the systems of SDEs (\ref{eqn:noncollBM})
and (\ref{eqn:noncollBESQ}) from the state $N \delta_0$ \cite{KT04}.

\SSC{Relaxation Phenomenon in Process 3 \label{sec:relax}}

As an application of Theorem \ref{thm:3_process},
here we study non-equilibrium dynamics of Process 3.
For Processes 1 and 2, see \cite{KT09,KT10,KT11}.

For Process 3 we consider the following special initial configuration
\begin{equation}
\eta(\cdot)=\sum_{j=1}^N \delta_{w_j}(\cdot)
\quad \mbox{with} \quad 
w_j=\frac{2 \pi r}{N}(j-1), \quad
1 \leq j \leq N.
\label{eqn:ed1}
\end{equation}
It is an unlabeled configuration with equidistant spacing
on $[0, 2 \pi r)$.
In this case the entire function for Process 3 
given by (\ref{eqn:Phi}) becomes
\begin{eqnarray}
\Phi_{\eta}^{w_k}(z)
&=& \prod_{\substack{1 \leq \ell \leq N, \cr \ell \not=k}}
\frac{\sin(z/2r-(\ell-1)\pi/N)}{\sin((k-\ell)\pi/N)}
\nonumber\\
&=& \prod_{n=1}^{N-1}
\frac{\sin[\{z/2r-(k-2)\pi/N\}+(n-1)\pi/N]}{\sin(n \pi/N)}.
\nonumber
\end{eqnarray}
We use the product formulas
$$
\prod_{n=1}^{N-1} \sin \left( \frac{n \pi}{N} \right) 
=\frac{N}{2^{N-1}},
\qquad
\prod_{n=1}^{N} \sin \left[ x+\frac{(n-1)\pi}{N} \right]
=\frac{\sin(Nx)}{2^{N-1}},
$$
and obtain
$$
\Phi_{\eta}^{w_k}(z)
=\frac{1}{N} \frac{\sin[N\{z-2 \pi r(k-1)/N\}/2r]}
{\sin[\{z-2\pi r(k-1)/N\}/2r]}.
$$

It is easy to confirm the equality
$$
\frac{\sin(Nx)}{\sin x}
=\sum_{\substack{m \in \Z, \cr |\sigma_N(m)| \leq (N-1)/2}}
e^{2 i \sigma_N(m) x},
\quad N \in \N,
$$
where $\sigma_N$ is defined by (\ref{eqn:spin}). 
Then, by (\ref{eqn:cal1}), the martingale function is given by
\begin{eqnarray}
\cM_{\eta}^{w_k}(t,y) &=& 
\int_{\R} d \widetilde{y} \,
\frac{1}{N} 
\sum_{\substack{m \in \Z, \cr |\sigma_N(m)| \leq (N-1)/2}}
e^{2 i \sigma_N(m)\{y+i\widetilde{y}-2 \pi r (k-1)/N\}/2r}
p_{\rm BM}(t, \widetilde{y}|0)
\nonumber\\
&=& \frac{1}{N} \sum_{\substack{m \in \Z, \cr |\sigma_N(m)|\leq (N-1)/2}}
e^{\sigma_N(m)^2 t/2r^2-i \sigma_N(m) y/r+2i (k-1) \sigma_N(m) \pi/N},
\nonumber
\end{eqnarray}
where $p_{\rm BM}$ is given by (\ref{eqn:p}).

Now the correlation kernel is obtained by the formula (\ref{eqn:K1}),
$$
\mbK_{\eta}(s,x;t,y)
= \cG_{\eta}(s,x;t,y)-\1(s>t) p^r(s-t, x|y; N)
$$
with
\begin{eqnarray}
\cG_{\eta}(s,x;t,y)
&=& \sum_{k=1}^N p^r(s, x| w_k; N)
\cM_{\eta}^{w_k}(t,y)
\nonumber\\
&=& \frac{1}{N} \frac{1}{2 \pi r}
\sum_{k=1}^N 
\sum_{\ell \in \Z} 
\sum_{\substack{m \in \Z, \cr |\sigma_N(m)|\leq (N-1)/2}}
{\cal K}_{k,\ell,m}(s,x;t,y),
\nonumber
\end{eqnarray}
where we use the expression (\ref{eqn:pR2})
for $p^r$ and obtain
$$
{\cal K}_{k,\ell,m}(s,x;t,y)
= e^{\{\sigma_N(m)^2 t - \sigma_N(\ell)^2 s\}/2r^2
+i \{\sigma_N(\ell) x- \sigma_N(m) y\}/r+2i (k-1)\{\sigma_N(m)-\sigma_N(\ell)\} \pi/N}.
$$
By the equality
$$
\sum_{k=1}^N e^{2i(k-1)\{\sigma_N(m)- \sigma_N(\ell)\} \pi/N}
=N \sum_{k \in \Z} \1(\ell=m+k N),
$$
we obtain the following decomposition,
$$
\cG_{\eta}(s, x; t, y)=\sum_{k \in \Z} \cG_{\eta}^{(k)}(s, x; t, y)
$$
with
\begin{eqnarray}
\cG_{\eta}^{(k)}(s, x; t, y)
&=& \frac{1}{2 \pi r}
\sum_{\substack{m \in \Z, \cr |\sigma_N(m)|\leq (N-1)/2}}
e^{-\{\sigma_N(m+kN)^2-\sigma_N(m)^2\}s/2r^2}
\nonumber\\
&& \qquad \qquad \qquad \times
e^{\sigma_N(m)^2(t-s)/2r^2-i \{\sigma_N(m)y-\sigma_N(m+kN)x\}/r}.
\nonumber
\end{eqnarray}
Since $\sigma_N(m+kN)^2 > \sigma_N(m)^2$ if $m \in \Z, |\sigma_N(m)| \leq (N-1)/2$
and $k \not=0$, we see that
for $(s, x), (t, y) \in [0, \infty) \times [0, 2 \pi r)$
$$
\lim_{T \to \infty} \cG_{\eta}^{(k)}(s+T, x; t+T, y)
= \left\{ \begin{array}{ll}
\cG_{\rm eq}(t-s,y-x),
\quad &\mbox{if $k=0$}, \cr
0, \quad &\mbox{otherwise},
\end{array} \right.
$$
where 
$$
\cG_{\rm eq}(t, x)
= \frac{1}{2 \pi r} 
\sum_{\substack{m \in \Z, \cr |\sigma_N(m)|\leq (N-1)/2}}
e^{\sigma_N(m)^2 t/2r^2 -i \sigma_N(m) x/r}.
$$
In particular, when $s=t$ we have
$$
\cG_{\rm eq}(0, x)
= \frac{1}{2 \pi r} 
\sum_{\substack{m \in \Z, \cr |\sigma_N(m)|\leq (N-1)/2}}
e^{-i \sigma_N(m) x/r}
= \frac{1}{2 \pi r} 
\frac{\sin(N x/2r)}{\sin(x/2r)}.
$$
The results are summarized as follows.

\begin{prop}
\label{thm:relax}
Let $(\Xi, \P_{\eta})$ be the Process 3
started at the configuration (\ref{eqn:ed1}). 
It shows a relaxation
phenomenon to an equilibrium determinantal process
$(\Xi, \P_{\rm eq})$.
The correlation kernel of $(\Xi, \P_{\rm eq})$ is given as
\begin{eqnarray}
&& \mbK_{\rm eq}^r(t-s, y-x)
= \cG_{\rm eq}(t-s, y-x)- \1(s>t) p^r(s-t, x|y; N)
\nonumber\\
&& \qquad = \left\{ \begin{array}{ll}
\displaystyle{ \frac{1}{2\pi r}
\sum_{\substack{m \in \Z, \cr |\sigma_N(m)|\leq (N-1)/2}}
e^{\sigma_N(m)^2(t-s)/2r^2-i \sigma_N(m)(y-x)/r}},
& \mbox{if $s < t$}, \cr
& \cr
\displaystyle{
\frac{1}{2 \pi r} \frac{\sin[N(y-x)/2r]}{\sin[(y-x)/2r]}},
& \mbox{if $s=t$}, \cr
& \cr
\displaystyle{ -\frac{1}{2\pi r}
\sum_{\substack{m \in \Z, \cr |\sigma_N(m)| > (N-1)/2}}
e^{\sigma_N(m)^2(t-s)/2r^2-i \sigma_N(m) (y-x)/r}},
& \mbox{if $s > t$}
\end{array} \right.
\nonumber
\end{eqnarray}
for $(s,x), (t,y) \in [0, \infty) \times [0, 2 \pi r)$.
\end{prop}

Let $\mu_{\rm eq}$ be the determinantal point process
with the spatial correlation kernel 
$$
\bK_{\rm eq}^r(y-x)
= \frac{1}{2 \pi r} \frac{\sin[N(y-x)/2r]}{\sin[(y-x)/2r]},
\quad x, y \in [0, 2 \pi r).
$$
This is the CUE eigenvalue distribution in the random
matrix theory (see, for instance, Section 11.1 in \cite{Meh04}).
The above equilibrium process $(\Xi, \P_{\rm eq})$
is reversible with respect to $\mu_{\rm eq}$.
This equilibrium determinantal process
was discussed also in \cite{NF02} as a stationary process
in several non-equilibrium dynamics which are expressed
using quaternion determinants (Pfaffians) and
are different from the present determinantal process $(\Xi, \P_{\eta})$.
\vskip 0.3cm

\noindent{\bf Remark 2} \,
The particle-density is given by
$\rho=N/2 \pi r$.
As expected, we can see 
$$
\lim_{\substack{r \to \infty, N \to \infty, \cr \rho={\rm const.}}}
\mbK^r_{\rm eq}(t-s,y-x)=\bK_{\rm sin}(t-s,y-x),
$$
where $\bK_{\rm sin}$ is the {\it extended sine kernel} 
with density $\rho$ defined on $\R$,
$$
\bK_{\rm sin}(t-s,y-x)
= \left\{ \begin{array}{ll}
\displaystyle{ \int_0^{\rho} dv \,
e^{\pi^2 v^2(t-s)/2} \cos(\pi v(y-x))},
& \mbox{if $s<t$}, \cr
& \cr
\displaystyle{ \frac{\sin(\pi \rho(y-x))}{\pi(y-x)}},
& \mbox{if $s=t$},\cr
& \cr
\displaystyle{ -\int_{\rho}^{\infty} dv \,
e^{\pi^2 v^2(t-s)/2} \cos(\pi v(y-x))},
& \mbox{if $s>t$}. \cr
\end{array} \right.
$$
\vskip 0.3cm
\noindent
It is the correlation kernel for the determinantal process
describing the equilibrium dynamics of
Process 1 with an infinite number of particles \cite{KT10}.

\vskip 0.3cm
\noindent{\bf Remark 3} \,
Consider the case that all particles have the same
constant drift with drift coefficient $b$.
The SDEs given by (\ref{eqn:noncollBM_S1}) are just modified by
$
B_j(t) \to B_j(t)+ bt,
$
that is, each BM is replaced by drifted Brownian motion.
In the correlation kernel (\ref{eqn:K1})
for Process 3, 
what we have to do is only replacing
$p^r(t,y|x; N)$ by its drift transform
$
p^r_b(t, y|x; N)
\equiv p^r(t,y|x; N) e^{b(y-x)-b^2 t/2}
= p^r(t,y-bt|x; N)
$
and $\cM_{\xi}^v(t, y)$ by
$\cM_{\xi}^v(t, y-bt)$. 
The correlation kernel 
is then given by
\begin{eqnarray}
\mbK_{\xi,b}(s,x; t,y)
&=& \int_{[0, 2 \pi r)} \xi(dv) p^{r}_b(s,x|v; N)
\cM_{\xi}^v(t, y-bt)- \1(s>t) p^{r}_b(s-t, x|y; N)
\nonumber\\
&=& \mbK_{\xi}(s, x-bs; t, y-bt),
\nonumber
\end{eqnarray}
$(s,x), (t,y) \in [0, \infty) \times [0, 2 \pi r)$.
The correlation kernel in the equilibrium is obtained by
$
\mbK_{{\rm eq},b}^r(t-s, y-x)
=\mbK_{\rm eq}^r(t-s, y-x-b (t-s)).
$

\SSC{Complex-Process Representations for BES$^{(n+1/2)}$ \label{sec:CPR}}

Let $\widehat{Y}^{(\nu)}(t), t \in [0, \infty)$ be the
Bessel process with index $\nu >-1$.
It solves the SDE
$$
\widehat{Y}^{(\nu)}(t)
= u+ B(t) + \frac{2\nu+1}{2} \int_0^t \frac{ds}{\widehat{Y}^{(\nu)}(s)},
\quad u >0, 
\quad t \geq 0,
$$
where if $-1 < \nu < 0$ the reflection boundary condition
is assumed at the origin.
The TD
is obtained from (\ref{eqn:p_nu}) as
\begin{eqnarray}
\widehat{p}^{(\nu)}(t, y|x) &=& p^{(\nu)}(t, y^2|x^2) 2y
\nonumber\\
&=& \left\{ \begin{array}{ll}
\displaystyle{
\frac{y^{\nu+1}}{x^{\nu}}
\frac{e^{-(x^2+y^2)/2t}}{t}
I_{\nu} \left( \frac{xy}{t} \right)},
& \quad t>0, x>0, y \in \R_+, \cr
\displaystyle{
\frac{y^{2\nu+1}e^{-y^2/2t}}{2^{\nu} t^{\nu+1} \Gamma(\nu+1)}},
& \quad t >0, x=0, y \in \R_+, \cr
& \cr
\delta_{x}(\{y\}),
& \quad t=0, x, y \in \R_+.
\end{array} \right.
\nonumber
\end{eqnarray}
Corresponding to this,
the integral transformation (\ref{eqn:I}) 
for Process 2 associated with BESQ$^{(\nu)}$
is converted into that
for BES$^{(\nu)}$ as 
\begin{equation}
\widehat{\cI}^{(\nu)}[\widehat{f}(W)|(t,x)]
=\int_{\R_+} dw \, \widehat{f}(i w) \widehat{q}^{(\nu)}(t,w|x),
\quad (t, x) \in [0, \infty) \times \R_+,
\label{eqn:transBES1}
\end{equation}
where
\begin{equation}
\widehat{q}^{(\nu)}(t, w|x)
= \frac{w^{\nu+1}}{x^{\nu}}
\frac{e^{(x^2-w^2)/2t}}{t}
J_{\nu} \left( \frac{x w}{t} \right),
\quad t > 0, \quad x>0, 
\label{eqn:q4}
\end{equation}
(some arguments will be given below for $x=0$ and 
$t=0$)
and
$\widehat{f}(z)$ is assumed to be a polynomial of $z^2$;
and thus 
\begin{equation}
\widehat{f}(-z)=\widehat{f}(z), \quad z \in \C.
\label{eqn:symf}
\end{equation}

Assume $x>0$. 
For $m \in \N_0$, the Bessel functions have
the following expansions by the trigonometric functions,
\begin{eqnarray}
J_{2m+1/2}(x)
&=& (-1)^m \sqrt{\frac{2}{\pi x}}
\left[ \sin x \sum_{k=0}^m
\frac{(-1)^k (2m+2k)!}{(2k)!(2m-2k)!} (2x)^{-2k} 
\right.
\nonumber\\
&& \qquad \qquad \quad \left.
+ \cos x \sum_{k=0}^{m-1}
\frac{(-1)^k (2m+2k+1)!}{(2k+1)!(2m-2k-1)!} (2x)^{-(2k+1)} \right],
\nonumber\\
J_{2m+3/2}(x)
&=& (-1)^m \sqrt{\frac{2}{\pi x}}
\left[ - \cos x \sum_{k=0}^m
\frac{(-1)^k (2m+2k+1)!}{(2k)!(2m-2k+1)!} (2x)^{-2k} 
\right.
\nonumber\\
&& \qquad \qquad \quad \left.
+ \sin x \sum_{k=0}^{m}
\frac{(-1)^k (2m+2k+2)!}{(2k+1)!(2m-2k)!} (2x)^{-(2k+1)} \right].
\label{eqn:Jodd}
\end{eqnarray}
They are obtained from Eq.(4.6.11) in \cite{AAR99}.
For example, if we set $m=0$ in (\ref{eqn:Jodd}), we have
$$
J_{1/2}(x)=\sqrt{\frac{2}{\pi x}} \sin x, \quad
J_{3/2}(x)=\sqrt{\frac{2}{\pi x}} 
\left( \frac{\sin x}{x} - \cos x \right).
$$
Then (\ref{eqn:transBES1}) with (\ref{eqn:q4}) gives 
\begin{eqnarray}
\widehat{\cI}^{(1/2)} \left[ \left. \widehat{f}(W) 
\right| (t,x) \right]
&=& \int_{\R_+} dw \, \widehat{f}(iw) \sqrt{\frac{2}{\pi t}} \frac{w}{x}
e^{-(-x^2+w^2)/2t} \sin (xw/t)
\nonumber\\
&=& \frac{1}{\sqrt{2 \pi t}} \frac{1}{i x} 
\int_{\R_+} dw \, \widehat{f}(iw) w
\left\{ e^{-(w-ix)^2/2t}-e^{-(w+ix)^2/2t} \right\},
\nonumber
\end{eqnarray}
and
\begin{eqnarray}
&& 
\widehat{\cI}^{(3/2)} \left[ \left. \widehat{f}(W) 
\right| (t,x) \right]
\nonumber\\
&& \qquad 
= \int_{\R_+} dw \, \widehat{f}(iw) \sqrt{\frac{2}{\pi t}} \frac{w}{x}
e^{-(-x^2+w^2)/2t} 
\left\{ \frac{t}{xw} \sin (xw/t) - \cos (xw/t) \right\}
\nonumber\\
&& \qquad
= \frac{1}{\sqrt{2 \pi t}}
\left[ \frac{t}{i x^3} \int_{\R_+} dw \, \widehat{f}(iw) w
\left\{ e^{-(w-ix)^2/2t}-e^{-(w+ix)^2/2t} \right\}
\right.
\nonumber\\
&& \qquad \qquad \qquad \left.
- \frac{1}{x^2} 
\int_{\R_+} dw \, \widehat{f}(iw) w^2
\left\{ e^{-(w-ix)^2/2t}+e^{-(w+ix)^2/2t} \right\}
\right].
\nonumber
\end{eqnarray}
By the assumption (\ref{eqn:symf}), they are rewritten as
\begin{eqnarray}
\widehat{\cI}^{(1/2)} \left[ \left. \widehat{f}(W) 
\right| (t,x) \right]
&=&
\frac{1}{\sqrt{2 \pi t}}
\frac{1}{(-i)x} \int_{\R} dw \, \widehat{f}(iw) w
e^{-(w+ix)^2/2t}
\nonumber\\
&=& 
\int_{-\infty+ix}^{\infty+ix} d \widetilde{y} \,
\frac{x+i \widetilde{y}}{x} \widehat{f}(x+i\widetilde{y})
p_{\rm BM}(t,\widetilde{y}|0),
\nonumber\\
\widehat{\cI}^{(3/2)} \left[ \left. \widehat{f}(W) 
\right| (t,x) \right]
&=&
\frac{1}{\sqrt{2 \pi t}}
\left[ \frac{t}{(-i) x^3}
\int_{\R} dw \,  \widehat{f}(i w) w e^{-(w+ix)^2/2t}
\right.
\nonumber\\
&& \qquad \qquad  \left.
+ \frac{1}{(-i)^2 x^2}
\int_{\R} dw \, \widehat{f}(iw) w^2 e^{-(w+ix)^2/2t} \right]
\nonumber\\
&=& 
\int_{-\infty+ix}^{\infty+ix} d \widetilde{y} \,
\left\{ \frac{t(x+i\widetilde{y})}{x^3}+\frac{(x+i\widetilde{y})^2}{x^2} \right\}
\widehat{f}(x+i\widetilde{y})
p_{\rm BM}(t,\widetilde{y}|0),
\nonumber
\end{eqnarray}
where we have changed the integral variables by
$\widetilde{y}=w+ix$.
Since the integrands are entire,
$\int_{-\infty +ix}^{\infty +ix} du \, (\cdot)$ can be
replaced by $\int_{\R} du \, (\cdot)$ by Cauchy's integral theorem.
Then we have the following expressions for the integral transformations
(\ref{eqn:transBES1}) with $\nu=1/2$ and 3/2, 
\begin{eqnarray}
&& \widehat{\cI}^{(1/2)} \left[ \left. \widehat{f}(W) \right|
(t,x) \right]
=
\widetilde{\rE} \left[ 
\frac{x+i \widetilde{Y}}{x} 
\widehat{f}(x+i \widetilde{Y})
\right], 
\nonumber\\
&&
\widehat{\cI}^{(3/2)} \left[ \left. \widehat{f}(W) \right|
(t,x) \right]
=
\widetilde{\rE} \left[
\left\{ \frac{t (x+i \widetilde{Y})}{x^3}
+\frac{(x+i \widetilde{Y})^2}{x^2} \right\}
\widehat{f}(x+i \widetilde{Y}) \right], 
\nonumber
\end{eqnarray}
$(x,t) \in [0, \infty) \times \R_+$.
These calculations are generalized as follows.

\begin{lem}
\label{thm:BES_mart1}
Let $\widehat{f}(z)$ be a polynomial of $z^2$.
Then for $(t,x) \in [0, \infty) \times \R_+$
\begin{equation}
\widehat{\cI}^{(n+1/2)} \left[ \left. \widehat{f}(W) \right|
(t,x) \right]
=\widetilde{\rE} \left[
Q_t^{(n+1/2)}(x+i \widetilde{Y}(t)) 
\widehat{f}(x+i \widetilde{Y}(t)) \right],
\quad n \in \N_0,
\label{eqn:BES_mart1}
\end{equation}
where
\begin{equation}
Q^{(n+1/2)}_t(x+iy) =
\left\{ \begin{array}{ll}
\displaystyle{\left( \frac{t}{2} \right)^n
\frac{x+iy}{x^{2n+1}}
\sum_{k=0}^n \frac{(2n-k)!}{(n-k)! k!}
\left( \frac{2 x (x+iy)}{t} \right)^k}, 
& \quad x>0, y \in \R, \cr
\displaystyle{
\frac{\sqrt{\pi}}{2^{n+1} t^{n+1} \Gamma(n+3/2)}
({\rm sgn} \, y)^n y^{n+2},
}
& \quad x=0, y \in \R,
\end{array} \right.
\label{eqn:Q1}
\end{equation}
for $t>0$ and $Q_0^{(n+1/2)} \equiv 1$.
\end{lem}
\noindent{\it Proof.} \,
For $x>0, t>0$, 
the integral kernel $\widehat{q}^{(n+1/2)}, n \in \N_0$
given by (\ref{eqn:q4}) is expanded following (\ref{eqn:Jodd})
and the integral (\ref{eqn:transBES1}) is calculated.
As demonstrated above for $n=0$ and 1,
integrals of functions including $\sin(xw/t)$ and $\cos(xw/t)$
over $\R_+$ are rewritten by using the integrals
of functions including $e^{-(w+ix)^2/2t}$ over $\R$.
Then we change the integral variable from $w$ to $\widetilde{y}$
by $\widetilde{y}=w+ix$.
Since all the integrands are entire with respect to $\widetilde{y}$,
integrals on a line $(-\infty+ix, \infty+ix)$ in $\C$ 
can be replaced by those on the real axis $\R$.
Then we have the integrals in the form
$\int_{\R} d\widetilde{y} 
Q_t^{(n+1/2)}(x+i\widetilde{y}) \widehat{f}(x+i\widetilde{y})
p_{\rm BM}(t,\widetilde{y}|0)$
with the first line of (\ref{eqn:Q1}), 
which is expressed as RHS of (\ref{eqn:BES_mart1})
by definition of $\widetilde{\rE}$.
The function $J_{\nu}(z)/z^{\nu}$ is entire and
(\ref{eqn:q4}) gives
$$
\widehat{q}^{(\nu)}(t,w|0)=\lim_{x \downarrow 0} \widehat{q}^{(\nu)}(t,w|x)
=\frac{w^{2\nu+1}e^{-w^2/2t}}{2^{\nu} t^{\nu+1} \Gamma(\nu+1)},
\quad t > 0.
$$
Then we obtain
$$
\widehat{\cI}^{(n+1/2)}[ \widehat{f}(W)|(t,0)]
=\frac{\sqrt{\pi}}{2^{n+1} t^{n+1} \Gamma(n+3/2)}
\int_{\R} d\widetilde{y} \,
({\rm \, sgn} \, \widetilde{y})^n \widetilde{y}^{n+2} \widehat{f}(i\widetilde{y})
p_{\rm BM}(t,\widetilde{y}|0),
\quad t >0,
$$
which gives the second line of (\ref{eqn:Q1}). 
When $t=0$, (\ref{eqn:BES_mart1}) gives $\widehat{f}(x)$ with $Q_0^{(n+1/2)}\equiv 1$.
Then the proof is completed.\qed
\vskip 0.3cm
When $\nu=-1/2$, $J_{-1/2}(x)=\sqrt{2/\pi x} \cos x$,
and we find
\begin{eqnarray}
\widehat{\cI}^{(-1/2)} \left[ \left. \widehat{f}(W) \right|
(t,x) \right]
&=& \frac{1}{\sqrt{2 \pi t}}
\int_{\R_+} d w \, \widehat{f}(i w)
\{e^{-(w-ix)^2/2t} + e^{-(w+ix)^2/2t} \}
\nonumber\\
&=& \int_{\R} d \widetilde{y} \,
\widehat{f}(x+i \widetilde{y}) p_{\rm BM}(t, \widetilde{y}|0) 
= \widetilde{\rE}[\widehat{f}(x+i \widetilde{Y}(t))].
\label{eqn:nu-1/2}
\end{eqnarray}
Then we define 
\begin{equation}
Q^{(-1/2)}_t \equiv 1.
\label{eqn:Q-1/2}
\end{equation}

Now we consider the following interacting particle system.
\begin{description}
\item[Process 2$^{\prime}$:] \,
Noncolliding Bessel process (BES$^{(\nu)}$) with $\nu >-1$;
\begin{eqnarray}
S &=& \R_+ \equiv \{x \in \R: x \geq 0\},
\nonumber\\
\widehat{\Xi}^{(\nu)}(t, \cdot)
&=& \sum_{j=1}^N \delta_{\widehat{X}^{(\nu)}_j(t)} (\cdot),
\quad t \in [0, \infty),
\nonumber\\
\widehat{X}^{(\nu)}_j(t)
&=& u_j+ B_j(t)
+ \frac{2 \nu+1}{2} \int_0^t \frac{ds}{\widehat{X}^{(\nu)}_j(s)}
\nonumber\\
&& 
+ \sum_{1 \leq k \leq N, k \not=j} 
\left\{ \int_0^t \frac{ds}{\widehat{X}^{(\nu)}_j(s)-\widehat{X}^{(\nu)}_k(s)}
+\int_0^t \frac{ds}{\widehat{X}^{(\nu)}_j(s)+\widehat{X}^{(\nu)}_k(s)} \right\},
\nonumber
\end{eqnarray}
$1 \leq j \leq N, t \geq 0$,
where $\xi=\sum_{j=1}^N \delta_{u_j} \in \mM_0(\R_+)$, 
$B_j(t), 1 \leq j \leq N, t \geq 0$
are independent BMs started at 0, 
and, 
if $-1 < \nu < 0$,
the reflection boundary condition is assumed at
the origin.
\end{description}
The process is denoted by
 $(\widehat{\Xi}^{(\nu)}, \widehat{\P}_{\xi}^{(\nu)}), \nu >-1$,
with an initial configuration $\xi \in \mM_0(\R_+)$
and expectation with respect to $\widehat{\P}_{\xi}^{(\nu)}$
is written as $\widehat{\E}^{(\nu)}_{\xi}$.

Since the transform from Process 2 to Process 2$^{\prime}$ is obvious,
as a corollary of Theorem \ref{thm:3_process} for Process 2,
CPR for the martingales
follows Lemma \ref{thm:BES_mart1} and (\ref{eqn:nu-1/2}) with (\ref{eqn:Q-1/2}) as, 
\begin{equation}
\widehat{\cM}_{\xi}^{u_k}(t, \widehat{Y}^{(n+1/2)}(t))
=\widetilde{\rE} \Big[ Q_t^{(n+1/2)}(\widehat{Z}^{(n+1/2)}(t))
\widehat{\Phi}_{\xi}^{u_k}(\widehat{Z}^{(n+1/2)}(t)) \Big],
\quad 1 \leq k \leq N,
\label{eqn:hatM}
\end{equation}
when $n \in \N_0 \cup \{-1\}$.
Here the complex process $\widehat{Z}^{(n+1/2)}(\cdot)$
is defined by (\ref{eqn:hatZ}), 
$\xi=\sum_{k=1}^N \delta_{u_k} \in \mM_0(\R_+)$, and
\begin{equation}
\widehat{\Phi}_{\xi}^{u_k}(z)
=\prod_{\substack{1 \leq \ell \leq N, \cr \ell \not=k}}
\frac{z^2-u_{\ell}^2}{u_k^2-u_{\ell}^2},
\quad z \in \C, \quad 1 \leq k \leq N.
\label{eqn:Phi2z}
\end{equation}
Note that BES$^{(-1/2)}$, $\widehat{Y}^{(-1/2)}(\cdot)$ is the reflecting 
Brownian motion on $\R_+$ with the reflecting wall at the origin.

Now CPR is given for Process 2$^{\prime}$ as follows.
Let $\widehat{Z}_j^{(n+1/2)}(t),  t \in [0, \infty), 1 \leq j \leq N$
be independent copies of (\ref{eqn:hatZ})
and the expectation with respect to them is
written as $\widehat{\bE}_{\u}$
when $\widehat{Z}_j^{(n+1/2)}(0)=u_j, 1 \leq j \leq N$.
\begin{prop}
\label{thm:BES}
Suppose 
$\xi=\sum_{j=1}^{N} \delta_{u_j} 
\in \mM_0(\R_+)$,
$\u=(u_1, \dots, u_N) \in \W_N(\R_+)$. 
Let $n \in \N_0 \cup \{-1\}$ and
$0 \leq t \leq T < \infty$.
For any ${\cal F}_{\widehat{\Xi}}(t)$-measurable bounded function $F$
we have
\begin{eqnarray}
&& \widehat{\E}_{\xi} \left[ F 
\left(\widehat{\Xi}^{(n+1/2)}(\cdot) \right) \right]
\nonumber\\
&&
=\widehat{\bE}_{\u}
 \left[F \left( \sum_{j=1}^{N} \delta_{\Re \widehat{Z}^{(n+1/2)}_j(\cdot)} \right)
\det_{1 \leq j, k \leq N}
\Big[
Q_T^{(n+1/2)}(\widehat{Z}^{(n+1/2)}_j(T))
\widehat{\Phi}_{\xi}^{u_k}
(\widehat{Z}^{(n+1/2)}_j(T)) \Big]
\right],
\nonumber
\end{eqnarray}
where $Q_t^{(n+1/2)}$ and 
$\widehat{\Phi}_{\xi}=\{\widehat{\Phi}_{\xi}^{u_k}, 1 \leq k \leq N \}$ are given by
(\ref{eqn:Q1}) with (\ref{eqn:Q-1/2}) and (\ref{eqn:Phi2z}), respectively.
\end{prop}
\vskip 0.3cm

In the present paper, we have focused on the connection between the
$h$-transforms and determinantal correlation functions
through DMR.
CPR obtained from DMR
as shown above will be, however, itself useful
to analyze the process.
In \cite{KT13}, the CBMR for Process 1
was used to prove the tightness of a series of processes,
which converges to the Dyson model with an infinite number of particles,
and the noncolliding property of the limit process.
Applications of CPRs for BES$^{(n+1/2)}, n \in \N_0$
will be reported elsewhere.

\vskip 0.8cm
\noindent{\bf Acknowledgements} \quad
The present author would like to thank 
T. Shirai, H. Osada, H. Tanemura, S. Esaki,
and P. Graczyk
for useful discussions.
He also thanks an anonymous referee for the 
careful reading of the manuscript,
and for a number of useful comments and suggestions
which improved the paper.
This work is supported in part by
the Grant-in-Aid for Scientific Research (C)
(No.21540397 and No.26400405) of Japan Society for
the Promotion of Science.

\end{document}